%% file: main.tex
\documentclass{article}

\usepackage[letterpaper,margin=1in]{geometry}
\input{preamble}

\input{newMacros}

\usepackage{pgfplots}
\pgfplotsset{compat=newest}
\pgfplotsset{scaled y ticks=false}
\usepgfplotslibrary{groupplots}
\usepgfplotslibrary{dateplot}
\tikzstyle{every node}=[font=\small]
\pgfplotsset{
    yticklabel style={/pgf/number format/fixed},  
}

\pgfplotsset{compat=1.11,
 /pgfplots/ybar legend/.style={
 /pgfplots/legend image code/.code={
 \draw[##1,/tikz/.cd,yshift=-0.25em]
 (0cm,0cm) rectangle (3pt,0.8em);},
 },
}

\title{Sequential change detection via backward confidence sequences}
\date{}
\usepackage{authblk}
\author[1]{Shubhanshu Shekhar \thanks{shubhan2@andrew.cmu.edu}}
\author[1,2]{Aaditya Ramdas \thanks{aramdas@stat.cmu.edu}}

\affil[1]{Department of Statistics and Data Science, Carnegie Mellon University}
\affil[2]{Machine Learning Department, Carnegie Mellon University}

\begin{document}
\maketitle
\input{abstract}
\tableofcontents

\input{main_body}

\end{document}

%% file: preamble.tex
\usepackage{amsmath,amssymb,amsthm,amsfonts,latexsym,bbm,xspace,graphicx,float,mathtools,mathdots}
\usepackage{braket,caption,subcaption,ellipsis,xcolor,textcomp}
\usepackage{combelow} %

\usepackage{hyperref}
\usepackage[nameinlink]{cleveref}
\crefname{ineq}{inequality}{inequalities}
\creflabelformat{ineq}{#2{\upshape(#1)}#3}

\usepackage{enumitem} 

\newtheorem{theorem}{Theorem}

\newtheorem*{namedtheorem}{\theoremname}
\newcommand{\theoremname}{testing}

\newtheorem{lemma}[theorem]{Lemma}

\newtheorem{proposition}[theorem]{Proposition}

\newtheorem{corollary}[theorem]{Corollary}

\theoremstyle{definition}
\newtheorem{definition}[theorem]{Definition}
\newtheorem{remark}[theorem]{Remark}

\newtheorem{example}[theorem]{Example}

\renewcommand{\Pr}{\mathop{\mathbb{P}\/}}

\newcommand{\bone}{\boldsymbol{1}}

\newcommand{\ignore}[1]{}

\usepackage{autonum}

\usepackage{natbib}
\setcitestyle{authoryear,open={(},close={)}} %

\newcommand{\NN}{\mathbb{N}}
\newcommand{\EE}{\mathbb{E}}

\DeclareMathOperator*{\argmax}{argmax}
\DeclareMathOperator*{\esssup}{esssup}

%% file: newMacros.tex
\newcommand{\iid}{\text{i.i.d.}\xspace}
\newcommand{\prob}{\mathbb{P}}
\newcommand{\expec}{\mathbb{E}}

\newcommand{\mc}[1]{\mathcal{#1}}
\newcommand{\mbb}[1]{\mathbb{#1}}
\newcommand{\lb}{\left[}
\newcommand{\rb}{\right]}
\newcommand{\lp}{\left(}
\newcommand{\rp}{\right)}

\newcommand{\defined}{\coloneqq}

\newcommand{\arl}{\text{ARL}\xspace} %
\newcommand{\hatT}{\widehat{T}} %
\newcommand{\hatepsilon}{\widehat{\epsilon}} %
\newcommand{\barB}{\bar{B}} %
\newcommand{\cusum}{\text{CuSum}\xspace}

\newcommand{\backwardtest}{\tau^{\text{back}}}
\newcommand{\dkl}{d_{\text{KL}}}
\newcommand{\dks}{d_{\text{KS}}}
\newcommand{\dmmd}{d_{\text{MMD}}}

\newcommand{\thetahat}{\widehat{\theta}}
\newcommand{\muhat}{\widehat{\mu}}
\newcommand{\sigmahat}{\widehat{\sigma}}
\newcommand{\sigmatilde}{\widetilde{\sigma}}
\newcommand{\Phat}{\widehat{P}}
\newcommand{\Qhat}{\widehat{Q}}
\newcommand{\fcsdetector}{\texttt{FCS-Detector}\xspace}
\newcommand{\bcsdetector}{\texttt{BCS-Detector}\xspace}

%% file: abstract.tex
\begin{abstract}
We present a simple reduction from sequential estimation to sequential changepoint detection (SCD). In short, suppose we are interested in detecting changepoints in some parameter or functional $\theta$ of the underlying distribution. We demonstrate that if we can construct a confidence sequence (CS) for $\theta$, then we can also successfully perform SCD for $\theta$. This is accomplished by checking if two CSs --- one forwards and the other backwards --- ever fail to intersect.  Since the literature on CSs has been rapidly evolving recently, the reduction provided in this paper immediately solves several old and new change detection problems. Further, our ``backward CS'', constructed by reversing time, is new and potentially of independent interest. We provide strong nonasymptotic guarantees on the frequency of false alarms and detection delay, and demonstrate numerical effectiveness on several problems.
\end{abstract}

%% file: main_body.tex
\section{Introduction}
\label{sec:introduction}

    We study the problem of sequential changepoint detection~(SCD), where the goal is to quickly detect any changes in the distribution generating a stream of observations, while controlling the false alarm rate at a specified level. Formally, for some (possibly infinite-dimensional) index set $\Theta$,  let $\{P_{\theta}\}_{\theta \in \Theta}$ denote a class of distributions  on some observation space $\mathcal{X}$. Suppose that for some $T \geq 1$, the observations $\{X_t: 1 \leq t \leq T\}$ are drawn \iid from $P_{\theta_0}$ with  $\theta_0\in \Theta$, and $\{X_t: t > T\}$ are drawn \iid from $P_{\theta_1}$ for some $\theta_1 \in \Theta$, with $\theta_1\neq \theta_0$.  Then, the SCD problem involves deciding between the null $H_0: \{T = \infty\}$, meaning no change occurred, and the alternative $H_1 = \cup_{i \in \mbb{N}} \{T = i\}$.

    Since the observations arrive sequentially, our task is to design a random stopping time, $\tau$, adapted to the natural filtration $\{\mc{F}_t: t \geq 1\}$ with $\mc{F}_t = \sigma(X_1, \ldots, X_t)$, at which we reject the null. A good stopping rule $\tau$ takes large values under the  the null~(i.e., when $T=\infty$), while minimizing the time required to detect the change under the alternative~(i.e., when $T<\infty$).
    Formally, when $T=\infty$, we require the \emph{average run length~(ARL)}, $\expec_\infty[\tau]$, to be  lower bounded by $1/\alpha$, for a given $\alpha \in (0,1]$, while ensuring a small \emph{detection delay}, $\mathbb{E}_T[(\tau-T)^+]$, when $T<\infty$. Informally, this means that we will have a false alarm roughly every $1/\alpha$ steps, so the reader may use $\alpha=10^{-3}$ as a rough guideline. (We also briefly discuss how to keep the probability of even a single false alarm below $\alpha$, but there is a tradeoff between the false alarm guarantee and  detection delay; detecting true changes quickly necessitates tolerating infrequent false alarms.)
    
    The literature on the topic of sequential changepoint detection is vast, as this problem arises in several important real-world applications, such as quality control~\citep{shewhart1930economic}, monitoring power networks~\citep{chen2015quickest}, analysis of genomes~\citep{chen2011bayesian, shen2012change}, and epidemic detection~\citep{baron2004early, yu2013change}. 
    Some of the earliest works in this topic \citep{shewhart1925application, page1954continuous, shiryaev1963optimum} assume that the pre- and post-change distributions admit known densities $f_0$ and $f_1$~(w.r.t. some common reference measure).
    These methods use statistics involving likelihood ratios, that can be computed efficiently in an incremental manner, and have also been shown to  admit strong optimality properties.
    The ideas underlying these likelihood-based schemes have also been extended to the case of (finite-dimensional) parametric families of distributions, such as the exponential family; see~\citet{tartakovsky2014sequential} for a detailed discussion. 
    However, these parametric assumptions are often too stringent to be applicable to many practical applications, where the data distributions may lie in much larger, nonparametric, classes. 
    Most of the ideas developed for the parametric setting, and in particular the likelihood-based schemes, fail to be applicable in the nonparametric case. With some exceptions discussed later, there are very few general principles for constructing nonparametric changepoint detection schemes. Our work in this paper addresses this issue, by developing a conceptually simple `meta-algorithm' for transforming any confidence sequence construction into a powerful changepoint detection method. As a consequence, we can immediately build upon the recent progress in constructing confidence sequences to instantiate new changepoint detection methods. 
    \begin{remark}
        \label{remark:same-class}
        We note that the SCD problem is usually studied in two settings, that differ from each other is a very subtle manner. In the  first~(and the more common) setting,  the pre- and post change distributions~($P_{\theta_0}$ and $P_{\theta_1}$) are assumed to lie in two different, and usually well-separated, classes of distributions. Using our notation, this is equivalent to assuming that there exist two known disjoint sets $\Theta_0$ and $\Theta_1$, such that $\theta_i \in \Theta_i$, for $i=0,1$.  The second setting, that is the subject of our paper, assumes less information. That is, both $\theta_0$ and $\theta_1$ are assumed to lie in some common index set~($\Theta$ in our notation), and the only condition is that $\theta_0 \neq \theta_1$. Hence, the first setting is in some sense ``easier'',  as the additional knowledge about $\Theta_0$ and $\Theta_1$~(their size, and their separation) can be exploited to design appropriate SCD schemes. While there exist some works that develop methods for SCD in the second setting, those schemes often rely on the specific structure of the problems considered. In this paper, we address this issue by developing a general principle for designing SCD schemes in the second setting. 
    \end{remark}
\section{Preliminaries}
\label{sec:preliminaries}
    The primary technical tool we use in our strategy are time-uniform version of confidence sets, called \emph{confidence sequences}~(CSs), that were first introduced in statistics literature by~\citet{darling1967confidence}. We present a definition adapted to our problem below. 
    \begin{definition}[\texttt{Confidence Sequences}]
        \label{def:confidence-sequences}
        Suppose $X_1, X_2, \ldots$ are drawn \iid from $P_\theta$, for some $\theta \in \Theta$. Then, for any $\alpha \in (0,1)$, a level-$(1-\alpha)$ CS, denoted by $\{C_t: t \geq 1\}$, is a collection of subsets $C_t \subset \Theta$, such that \textbf{(i)} $C_t$ is $\sigma(X_1,\dots,X_t)$-measurable and \textbf{(ii)}
            $\mathbb{P}\big( \forall t \geq 1: \theta \in C_t \big) \geq 1-\alpha$. 
    \end{definition}
    \begin{remark}
        \label{remark:nestedness} 
        Due to the time-uniformity in the definition of CSs, we can replace the confidence set $C_t$ with the smaller set $\widetilde{C}_t \defined \cap_{s\leq t} C_s$. The new CS, $\{\widetilde{C}_t: t \geq 1\}$, consists of \emph{nested} confidence sets; that is, $\widetilde{C}_t \subset \widetilde{C}_s$ for $s<t$. 
    \end{remark}
    \begin{remark}
       \label{remark:CS-def-general}
        The data do not need to be \iid for defining CSs. The above definition can be easily generalized to the case of independent random variables, with $X_t \sim P_{\vartheta_t}$, with $\vartheta_t \in \Theta$~(see~\Cref{appendix:background}). 
        This, however, requires that $\Theta$ is endowed with the notions of addition and scalar multiplication~(a sufficient condition is that $\Theta$ is a vector space), which we implicitly assume when needed. 
    \end{remark}
    We now introduce a notion of the `size' of the confidence set $C_t$, that reflects the amount of uncertainty. 
    \begin{definition}[\texttt{CS width}]
        \label{def:CS-width}
        Let $\Theta$ be endowed with a distance metric $d$, and let $\{C_t: t \geq 1\}$ denote a level-$(1-\alpha)$ CS constructed on observations $X_1, X_2, \ldots$ drawn \iid from $P_{\theta}$. A function $w(t, \theta, \alpha)$ denotes the pointwise width (bound) of the CS, if for all $t \in \mathbb{N}$ and $\theta \in \Theta$, we have $\sup_{\theta' \in C_t}d(\theta, \theta') \leq w(t, \theta, \alpha)$. We define the uniform width over $\Theta$, as $w(t, \Theta, \alpha) = \sup_{\theta \in \Theta} w(t, \theta, \alpha)$. 
    \end{definition}
    As we will see later in~\Cref{sec:instantiations}, most of the non-trivial CSs have their pointwise widths~(and often, the uniform widths as well) converging to $0$ with the number of observations. 
    \begin{example}
    \label{example:gaussian-mean-cs}
        Consider independent random variables $\{X_t: t \geq 1\}$, with $X_t \sim N(\theta, 1)$ and $\theta \in \mathbb{R}$ for all $t \geq 1$. In this case, the parameter set is $\Theta = \mathbb{R}$. For this process, we can define the CS $\{C_t: t \geq 1\}$ as follows: $C_t = [\bar{X}_t - w_t/2, \bar{X}_t + w_t/2]$, where $\bar{X}_t = (1/t) \sum_{i=1}^t X_i$, and $w_t = 3.4 \sqrt{  \lp \log \log (2t) + 0.72 \log(10.4/\alpha) \rp/t }$. 
        Thus, if we endow the parameter space $\Theta$ with the metric $d(\theta, \theta') = |\theta -\theta'|$, we observe that the uniform width of $C_t$, denoted by $w(t, \Theta, \alpha) = w_t$,  converges to $0$. 
    \end{example}
    
    \noindent\textbf{Related Work.}
        As we mentioned earlier, a large part of the existing SCD literature  focuses on the parametric setting. We refer the reader to some recent surveys, such as those by~\citet{veeravalli2014quickest, xie2021sequential}, and the textbook by~\citet{tartakovsky2014sequential} for  details. In this section, we discuss some results on nonparametric SCD methods that are more relevant to our work. 
    
        \citet{shin2022detectors} developed a novel framework for changepoint detection by introducing \emph{e-detectors}; obtained by combining a sequence of e-processes defined uniformly over the class of pre-change distribution. They showed that their resulting strategy using e-detectors controls the ARL under very general conditions, and also proved the optimality of their schemes (in terms of worst-case detection delays) in some cases. However, as we mentioned in~\Cref{remark:same-class}, their framework is applicable mainly in cases where the pre- and post-change distributions are known to lie in different classes. Our techniques, described in~\Cref{sec:proposed-scheme-1} and~\Cref{sec:proposed-scheme-2}, addresses this issue. 
        
        \citet{maillard2019sequential} considered the task of detecting a change in the mean of a sequence of independent, univariate, sub-Gaussian random variables; and proposed an SCD method by deriving a new, doubly time-uniform confidence sequence for the scan statistics associated with a generalized likelihood ratio scheme~\citep{lai2010sequential}. The original proof of this concentration inequality~\citep[Theorem 4]{maillard2019sequential} was incomplete, and a corrected version (with an additional $\log \log$ term)  was obtained by the author in~\citep[Chapter 3, \S~4.1]{maillard2019mathematics}. For the resulting scheme, \citet{maillard2019mathematics} obtained bounds on the probability of false alarm and on the detection delay, and also established the optimality of this scheme under certain scenarios~(such as for Gaussian observations). Unlike \citet{maillard2019mathematics}, our SCD framework is applicable to a much wider class of problems beyond univariate mean testing. Nevertheless, when specialized to the case of univariate sub-Gaussian observations, our scheme achieves the same detection delay bound, while also providing control over the ARL (instead of the probability of false alarm). 
        
        \citet{puchkin2022contrastive} considered the SCD problem under the assumption that both the pre- and post-change distributions admit densities w.r.t.\ a common reference measure, and proposed a strategy based on learning a discriminator to estimate the density ratio. They showed that their strategy can control ARL at the required level, and also obtained high probability upper bound on the detection delay in terms of the Jensen Shannon~(JS) divergence and the $L_2$ norm of the difference of densities. Unlike them, our framework does not require the existence of densities, and it also works for distributions that are separated in terms of a large class of metrics, and not just the JS divergence. 
        
        Another class of nonparametric schemes for SCD are based on the kernel-MMD metric, first employed by~\citet{gretton2012kernel} for designing powerful nonparametric two-sample tests. \citet{li2019scan} proposed a SCD scheme based on a variant of the block-MMD statistic~\citep{zaremba2013b} computed using the observations, and a block of pre-change data. More recently, \citet{flynn2019change} and~\citet{wei2022online} proposed new SCD schemes that use linear and block-MMD statistics to define nonparametric analogues of the CuSum test of~\citet{page1954continuous}. 
        However, these schemes suffer from weak theoretical guarantees on the detection delay. Furthermore, similar to the case of~\citet{puchkin2022contrastive}, the strategies for designing these SCD schemes are specific to the kernel-MMD metric; and there is no obvious way to extend them to other popular metrics, such as the Kolmogorov-Smirnov metric. Our work addresses these issues.

        Similarly, most other existing works in SCD are geared towards specific problem settings. Hence, both the design of the scheme as well as their analysis are strongly tied to the details of the problem being studied. Examples include empirical likelihood based methods for distributions on finite alphabets~\citep{lau2018binning}, nearest-neighbor techniques for multivariate or non-euclidean data~\citep{chen2019sequential}, and spectral scan statistics for graph valued data~\citep{sharpnack2013changepoint}. 
        However, our objective in this paper is different: instead of developing a powerful SCD scheme for a specific task, we  develop an abstract  unifying template for designing SCD schemes, that can then be instantiated for a large range of (old and new) SCD problems. 
        
   \noindent\textbf{Our Contributions.}
            In~\Cref{sec:proposed-scheme-1}, we first present (as a warmup) a changepoint detection scheme that uses a single level-$\alpha$ forward CS. Our strategy is to stop as soon as the CS becomes `inconsistent'; that is, it includes a point that it had previously discarded. We show in~\Cref{prop:general-strategy-1}, that this simple strategy controls the probability of false alarm at level $\alpha$, and we also obtain a high probability upper bound on its detection delay. 
            However, this scheme is too conservative as its ARL is infinite, and in practice this might result in large detection delays, especially when $T$ is large. 
            
            In~\Cref{sec:proposed-scheme-2}, we present our main strategy that proceeds by checking at each time $t$, whether a forward CS and a backward CS (a new notion) are consistent, and stops whenever an inconsistency is detected. 
            In~\Cref{theorem:general-strategy}, we show that the ARL of this scheme is at least $1/(2\alpha)$, and we  characterize its expected detection delay under general conditions. 
            
            Finally, in~\Cref{sec:instantiations}, we demonstrate the power and generality of our proposed scheme by instantiating it with five different confidence sequences. The general bound on the detection delay obtained in~\Cref{theorem:general-strategy} easily translate into problem-specific upper bounds in all these cases, and we also empirically verify the theoretical predictions through some simple numerical simulations. 
\section{\smash{Warmup: change detection via a forward CS}}
\label{sec:proposed-scheme-1}

    Before presenting our general scheme in the next section, we first introduce a simpler SCD method that only uses a single forward CS. We refer to this scheme as the \fcsdetector. The idea underlying this scheme is that if there is a change in the distribution generating the observations $\{X_t: t \geq 1\}$, then the intersection of the CS will eventually end up being empty. Formally,  we proceed by constructing a level-$(1-\alpha)$ confidence sequence~(CS) for the unknown $\theta_0$, denoted by $\{C_t: t \geq 1\}$, as introduced in~\Cref{def:confidence-sequences}.   When there is no changepoint, then the CS satisfies 
        $\prob_{\infty}(\forall t \in \NN: \theta_0 \in C_t) \geq 1-\alpha$.
    However, if there is a changepoint at some time $T$,  we expect that the confidence sets, $C_t$, deviate away from the confidence set $C_T$, for $t > T$. Eventually, after sufficiently many post-change observations, the confidence sequence will be inconsistent and self-contradictory. That is, at some time $t$ such that $t- T$ is large enough, we expect that $\cap_{s=1}^t C_s = \emptyset$. We thus define the stopping time, $\tau$, as the smallest $t$ at which the above inconsistency is observed. 
    \begin{definition}[\fcsdetector]
        \label{def:general-strategy-1}
        Given observations $X_1, X_2, \ldots$, we construct a confidence sequence~(CS), denoted by $\{C_t: t \geq 1\}$ for the pre-change parameter $\theta_0$. We stop at time 
            $\tau \defined \min \{n \geq 1: \exists t < n, \, C_t \cap C_n = \emptyset \}$. 
    \end{definition}
    This strategy  satisfies the following properties. 
    \begin{proposition}
    \label{prop:general-strategy-1}
         Consider a change point detection problem with observations $X_1, X_2, \ldots$ drawn \iid from $P_{\theta_0}$ for $t \leq T$ and from $P_{\theta_1}$ for $t>T$, with $T$ lying in $\mathbb{N} \cup \{\infty\}$ and $\theta_0, \theta_1 \in \Theta$. Suppose for any $\theta \in \Theta$, we can construct confidence sequences $\{C_t: t \geq 1\}$, with uniform width $ w(t) \equiv w(\cdot, \Theta, \alpha)$. 
         Then, we have the following: \\
            \textbf{(i)} When $T=\infty$, the \fcsdetector controls the probability of false alarm~(PFA) at level $\alpha$. That is, $\prob_{\infty} \lp \tau < \infty \rp \leq \alpha$.  \\
            \textbf{(ii)} Suppose $T$ is large enough to ensure that $w(T)< d(\theta_0, \theta_1)$. For $t >T$, define $\lambda_t := T/t$, and $\widetilde{\theta}_t := \lambda_t \theta_0 + (1-\lambda_t) \theta_1$. Then, we have 
                $\tau - T \leq   \min \left\{t-T:  w(t) + w(T) \leq d( \theta_0, \widetilde{\theta}_t ) \right\}$,
            with probability at least $1-\alpha$. 
    \end{proposition}
    
   \begin{remark}
        \label{remark:norm-induced}
        \sloppy If $d$  is induced by a norm on $\Theta$, we can  bound the detection delay under the event $\mc{E}$ by 
            $(\tau - T)^+ \leq \min \left\{ t - T: w(t) + w(T) \leq \frac{t-T}{t} d(\theta_0, \theta_1) \right\}$. 
        In many instances, we have $w(t) \approx 1/\sqrt{t}$ suppressing logarithmic factors. Then, the above expression implies that from `small' $T$, the delay is linear in $T$, while for large $T$, the delay behaves approximately like $\sqrt{T}$. The details of these calculations, along with plots of delay versus $T$ are in~\Cref{appendix:proof-general-strategy-1}. 
    \end{remark}
    \begin{remark}
        \label{remark:warmup-T-large-enough} 
        The condition on $T$ used above for obtaining the bound on detection delay is necessary, because we do not assume that the pre-change distribution~(i.e., the parameter $\theta_0$) is known to us. Thus, to be able to detect the change in distribution, we must have enough observations from the pre-change distribution to estimate $\theta_0$ accurately enough, in comparison the magnitude of change, $d(\theta_0, \theta_1)$. As an extreme example, if $T=1$, no method can realistically detect that a change occurred (since it is statistically plausible that all the data are simply i.i.d.\ and no change occurred at all). Said differently, $T=0$ and $T=\infty$ are information theoretically equivalent, and we need to be far enough away from those extremes for practical detectability.
    \end{remark}
    When $T=\infty$, the \fcsdetector continues sampling without stopping w.p.\ at least $1-\alpha$. Hence, its ARL  is infinite. This makes it is too conservative in detecting the changepoint --- we cannot provide an upper bound on the expected detection delay when $T<\infty$, and can only characterize the delay under an event of probability $1-\alpha$~(see~\Cref{fig:delay-strategy-1} in~\Cref{appendix:proof-general-strategy-1}). We address this next, by proposing a scheme that augments the forward CS with a series of backward CSs. 
\section{\smash{Change detection via a backward CS}}
\label{sec:proposed-scheme-2}
    We now introduce our main SCD strategy that addresses the two drawbacks of the simpler SCD strategy discussed in the previous section.  Informally, the idea underlying our strategy is as follows: in each round $t \geq 2$, we construct the usual forward CS, and a new backward CS~(using reversed observations, see~\Cref{def:backward-CS}). We refer to this scheme as the \bcsdetector. If there has been a changepoint, we expect the forward and backward CSs to concentrate on different regions of $\Theta$~(i.e., around $\theta_0$ and $\theta_1$ resp.). Hence, we stop as soon as they become inconsistent. See~\Cref{fig:illustrating-bcsdetector} in~\Cref{appendix:background} for a visual illustration, and~\Cref{appendix:repeated-sequential-tests} for an interpretation of our scheme in terms of repeated sequential tests.
    
    We now present our definition of  backward CSs. 
    \begin{definition}[Backward Confidence Sequences]
        \label{def:backward-CS}
        Let $X_1, X_2, \ldots$ be drawn \iid from $P_{\theta}$, for some $\theta \in \Theta$. For any $n \geq 1$, we say that a sequence of sets $\{B_t^{(n)}\}_{1\leq t \leq n} \subseteq \Theta$ is a backward CS, if  
            \textbf{(i)} $B_t^{(n)}$ is $\sigma\lp X_t, \ldots, X_n \rp$ measurable, and \textbf{(ii)} $\mathbb{P}\lp \forall t \in [n]: \theta \in B_t^{(n)} \rp \geq 1-\alpha$. 
    \end{definition}
    Note that for $n > 1$, a forward CS $C_t$ does not satisfy the first condition, since $C_t$ is built using $X_1,\dots,X_t$, but $B^{(n)}_t$ can only use $X_t, \dots, X_n$.
    But, a backward CS at any $n$ can be interpreted as the usual forward CS, introduced in~\Cref{def:confidence-sequences}, constructed on observations seen in a reverse order from $n$ to 1; see~\Cref{appendix:background} for details.

    Without loss of generality, we assume that any CS consists of a nested sequence of sets, as discussed in~\Cref{remark:nestedness}. 
    In other words, at a given time $n$,  $C_n$ is the smallest set among $\{C_t: t \in [n]\}$, while $B^{(n)}_1$ is the smallest among $\{B_t^{(n)}: t \in [n]\}$. 
    We say that `the two confidence sequences $\{C_t: t \geq 1\}$ and $\{B_t^{(n)}: t \in [n]\}$ `do not intersect at time $n$' if $C_r \cap B_s^{(n)} = \emptyset$ for some $1\leq r,s \leq n$. It is easy to check that  two nested CSs do not intersect if and only if their smallest sets (i.e., $C_n$ and $B_1^{(n)}$) do not intersect. 
    
    When $T=\infty$,  at every time $n$, both CSs $\{C_t: t \geq 1\}$ and $\{B_t^{(n)}: t \in [n]\}$ will  contain $\theta_0$ with probability at least $1-2\alpha$, and thus they will intersect with the same probability. 
    This motivates us to stop and declare a changepoint at the first time $n$ at which they do not intersect.
    \begin{definition}[\bcsdetector]
        \label{def:general-strategy}
        Given observations $X_1, X_2, \ldots$, suppose we construct a forward CS $\{C_t: t \geq 1\}$ and new backward CSs $\{B_t^{(n)}: t \in [n]\}$  for every $n \geq 1$. Assume that all the constructed CSs are nested. Then, we define the stopping time, $\tau$, as the first time at which the forward and backward CSs do not intersect:
            $\tau := \inf\{ n \geq 1: C_n \cap B_1^{(n)} = \emptyset\}$. 
    \end{definition}
    We now present the main result of this section. 
    \begin{theorem}
        \label{theorem:general-strategy}
        Consider a SCD problem with observations $X_1, X_2, \ldots$ drawn \iid from $P_{\theta_0}$ for $t \leq T$ and from $P_{\theta_1}$ for $t>T$, with $T$ lying in $\mathbb{N} \cup \{\infty\}$ and $\theta_0, \theta_1 \in \Theta$. Suppose for any $\theta \in \Theta$, we can construct confidence sequences $\{C_t: t \geq 1\}$ with pointwise width $w(\cdot, \theta, \alpha)$. 
        Then, we have the following: 
            \textbf{(i)} When, there is no changepoint ($T=\infty$), the \bcsdetector satisfies $\expec_{\infty}[\tau] \geq \frac{1}{2\alpha} - \frac{3}{2}$.  \\
            \textbf{(ii)} Suppose $T < \infty$, and the pre-change parameter $\theta_0$ is not known. Introduce the event $\mc{E} = \{\theta_0 \in C_t: 1 \leq t < T\}$, and note that $\prob_T\lp \mc{E} \rp \geq 1- \alpha$ by construction. Then, for $\alpha \in (0,0.5)$, we have $\expec_T\lb (\tau - T)^+ | \mc{E} \rb  \leq 3\frac{u_0(\theta_0, \theta_1, T)}{1-\alpha}$, where $u_0  \defined \min \{t \geq 1 :  w(t, \theta_1, \alpha) + w(T, \theta_0, \alpha)< d(\theta_1, \theta_0)\}$. 
    \end{theorem}
    Recall  from~\Cref{def:CS-width} that $w(t, \theta_1, \alpha)$  denotes the width of the level-$\alpha$ confidence set with $t$ observations, when the true parameter is $\theta_1$.
    The proof of this theorem is given in~\Cref{appendix:proof-general-strategy}. 
    In many problem instances, the pre-change parameter $\theta_0$ is known as it represents the `natural state' of the process being observed. We can specialize the above result stated  to this case as follows. 
    \begin{corollary}
        \label{corollary:known-pre-change}
        Suppose the pre-change parameter $\theta_0$ is known, and  $T<\infty$. Then, we have 
                $\expec_{T}[ \lp \tau - T \rp^+] \leq (3/(1-\alpha)) t_0(\theta_0, \theta_1, T)$,
            where  $t_0 \equiv t_0(\theta_0, \theta_1)  \defined \min \{t \geq 1 :  w(t, \theta_1, \alpha) < d(\theta_1, \theta_0)\}$. 
    \end{corollary}
    \begin{remark}
        \label{remark:comparison-of-the-two-schemes} These results  demonstrate how the drawbacks of the \fcsdetector~(introduced in~\Cref{sec:proposed-scheme-1}) are addressed by carefully incorporating the idea of backward CSs in the design strategy. In particular, this new scheme, called the \bcsdetector, is less conservative, and has a finite lower bound on the ARL under the null. More importantly, when $T<\infty$, the expected detection delay of \bcsdetector is also finite, and furthermore, it is also independent of the value of the changepoint $T$. This is in contrast to the (high probability) bound on the detection delay for~\fcsdetector, in which the detection delay increases approximately as $\sqrt{T}$, as the changepoint $T\to\infty$.  
    \end{remark}
    \begin{remark}[\texttt{Other estimates}]
        \label{remark:changepoint-estimate}
        We can also construct an estimate of the changepoint, denoted by $\hatT$, as the time at which $C_t$ and $B_t^{(\tau)}$ are most separated:
            $
            \hatT \equiv \hatT(\tau) \defined \max \; \argmax_{1 \leq t \leq \tau} \; d \big( C_t, B_t^{(\tau)} \big) = \max \{s \leq \tau: d(C_s, B_s^{(\tau)}) = \max_{1 \leq t\leq \tau} d(C_{t}, B_{t}^{(\tau)}) \}$. 
        Additionally, we define an estimate of the magnitude of the change $\epsilon \defined d(\theta_0, \theta_1)$ as the separation between $C_{\hatT}$ and $B_{\hatT}^{(\tau)}$, as measured by the distance metric $d$:
            $ 
            \hatepsilon \equiv \hatepsilon(\tau) \defined \max_{\theta \in C_{\hatT}, \theta' \in B_{\hatT}^{(\tau)}}\; d(\theta, \theta'). 
            $
        While we do not obtain theoretical guarantees, some empirical results in the next section indicate that these estimates are accurate for several instantiations of the \bcsdetector.  
    \end{remark}
\section{Instantiations of \bcsdetector}
\label{sec:instantiations}
    In the previous section, we introduced a conceptually simple device that allows us to transform any confidence sequence construction into a powerful, sequential changepoint detector. This allows us to instantiate our general changepoint detection meta-algorithm to various scenarios, by leveraging the recent progress in constructing confidence sequences. We illustrate this, by presenting a variety of parametric and nonparametric SCD problems in this section.

\subsection{Parametric Change of Mean Detection}
\label{subsec:univariate-gaussian}
    We begin by considering the simple case of univariate Gaussian mean changepoint detection~(although the same results also hold for the larger class of sub-Gaussian random variables). In this problem, we observe $\{X_t: t \geq 1\}$, drawn \iid according to the distribution $N(\mu_t, 1)$, with $\mu_t = \theta_0$ for $t \leq T$ and $\mu_t = \theta_1$ for $t>T$. Note that in this problem, we have $\Theta = \mathbb{R}$ and we can set the distance metric, $d$, to be the absolute value of the difference. We will use the CS for Gaussian means, recently derived by~\citet{howard2021time}, that we had introduced earlier in~\Cref{example:gaussian-mean-cs}. 
    Furthermore, we can use the same expression for constructing the backward CS at any time $n$, denoted by $\{B_t^{(n)}: t \geq 1\}$, but with the order of observations reversed, as described in~\Cref{sec:proposed-scheme-2}. The following result shows that in this parametric setting, our changepoint detection scheme achieves an order-optimal detection delay (i.e., optimal modulo poly-logarithmic factors): 
    \begin{corollary}
        \label{prop:gaussian-mean} 
        Suppose $\{X_t: t \geq 1\}$ are drawn \iid according to $P_{\theta_0} = N(\theta_0, 1)$ for $t\leq T$, and $P_{\theta_1} = N(\theta_1, 1)$ for $t > T$. Note that in this case, we have $\dkl(P_{\theta_1}, P_{\theta_0}) = (\theta_1-\theta_0)^2/2$. Then, if $T = \Omega \lp \log(1/\dkl(P_{\theta_1}, P_{\theta_0}))/ \dkl(P_{\theta_1}, P_{\theta_0})\rp$, then we have 
        \begin{align}
           {\small \mathbb{E}_T[(\tau-T)^+ | \mc{E}] = \mc{O}\lp \frac{ \log \log(\frac{1}{\dkl(P_{\theta_1}, P_{\theta_0})})  + \log(\frac{1}{\alpha})}{\dkl(P_{\theta_1}, P_{\theta_0})} \rp,   \label{eq:gaussian-mean-delay}}
        \end{align}
        where $\mc{E}$ is the $(1-\alpha)$ probability  event in~\Cref{theorem:general-strategy}. 
    \end{corollary}
    
    \begin{remark}
        \label{remark:gaussian-mean-near-optimality} 
        If the pre-change mean~($\theta_0$) is known, the above upper bound  holds for the worst-case detection delay without the conditioning,  defined as $J_L(\tau) \defined \sup_{T>0} \esssup \mathbb{E}[(\tau-T)^+|\mc{F}_{T}]$. Under the assumption that $\theta_1$ is also known,  \citet{lorden1971procedures} showed the following universal lower bound on this quantity: $\inf_{\tau'} J_L(\tau') = \frac{\log(1/\alpha)}{\dkl(P_{\theta_1}, P_{\theta_0})} \lp 1 + o(1) \rp$, as $\alpha \to 0$. Thus, \bcsdetector matches this optimal performance, modulo logarithmic factors, without the knowledge of the post-change parameter. 
        This is unlike some of the existing schemes for Gaussian mean change detection, such as~\citet{pollak1991sequential}, which achieve the same order optimal detection delay, but with additional assumptions~(known lower bound on change, and in the  limit of $T \to \infty$). 
    \end{remark}
    
    \noindent\textbf{Empirical Verification.} We now  verify the theoretical claims of our proposed changepoint detection scheme using observations  drawn from a unit-variance normal distribution with the pre-change mean $\theta_0=0$,  and post-change mean $\theta_1=\Delta$. 
    In~\Cref{fig:Gaussian-mean-CS}, we consider the case of $\Delta=0.4$ with the change occurring at $T=800$. The plots show the forward and backward CSs at the time at which the change is detected in a trial, as well as the distributions of the detection delay, estimated changepoint, and estimated change magnitude over $250$ trials. The plots indicate that \bcsdetector  detects changes quickly and accurately. 
    
   \iftrue 
    \begin{figure*}[tb]
            \def\figwidth{0.25\linewidth}
            \def\figheight{0.25\linewidth} %
            \centering
            \input{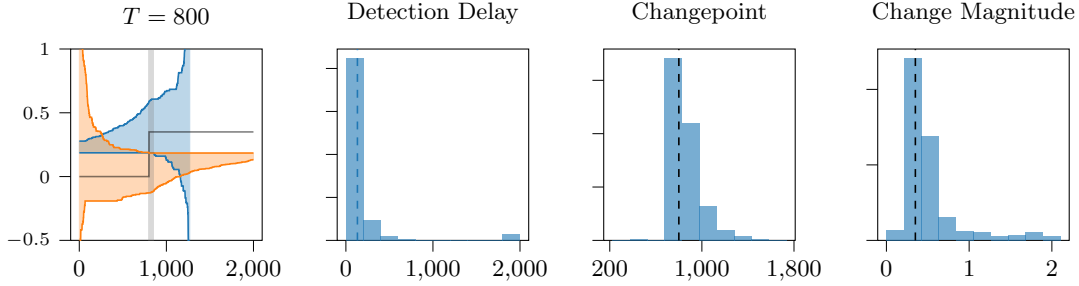} 
            \input{Figures/GaussianMean/GaussianMean5564DetectionDelay} 
            \input{Figures/GaussianMean/GaussianMean5564EstimatedChangepoint}  
            \input{Figures/GaussianMean/GaussianMean5564EstimatedChangeMagnitude}
            \caption{The figures show the performance of our changepoint detection scheme, \bcsdetector, with univariate Gaussian observations whose mean changes from $0$ to $0.4$ at the time $T=800$. The first plot shows the forward and backward CSs at time of detection~($\tau=1264$) in one of the trials, with the shaded gray region being the points at which the two CSs disagree. The next three plots show the empirical distribution of the detection delay, the estimated changepoint location, and the estimated changepoint magnitude over $250$ repeated trials.}
            \label{fig:Gaussian-mean-CS}
    \end{figure*}
   \fi

    In~\Cref{fig:Delay-vs-Delta}, we study the variation of the average detection delay  of our changepoint detection scheme as the change magnitude $\Delta$ is varied.  The empirical results verify the expected proportionality to $1/\Delta^2$ of the average detection delay, as claimed by our theoretical results. 
   \iftrue 
    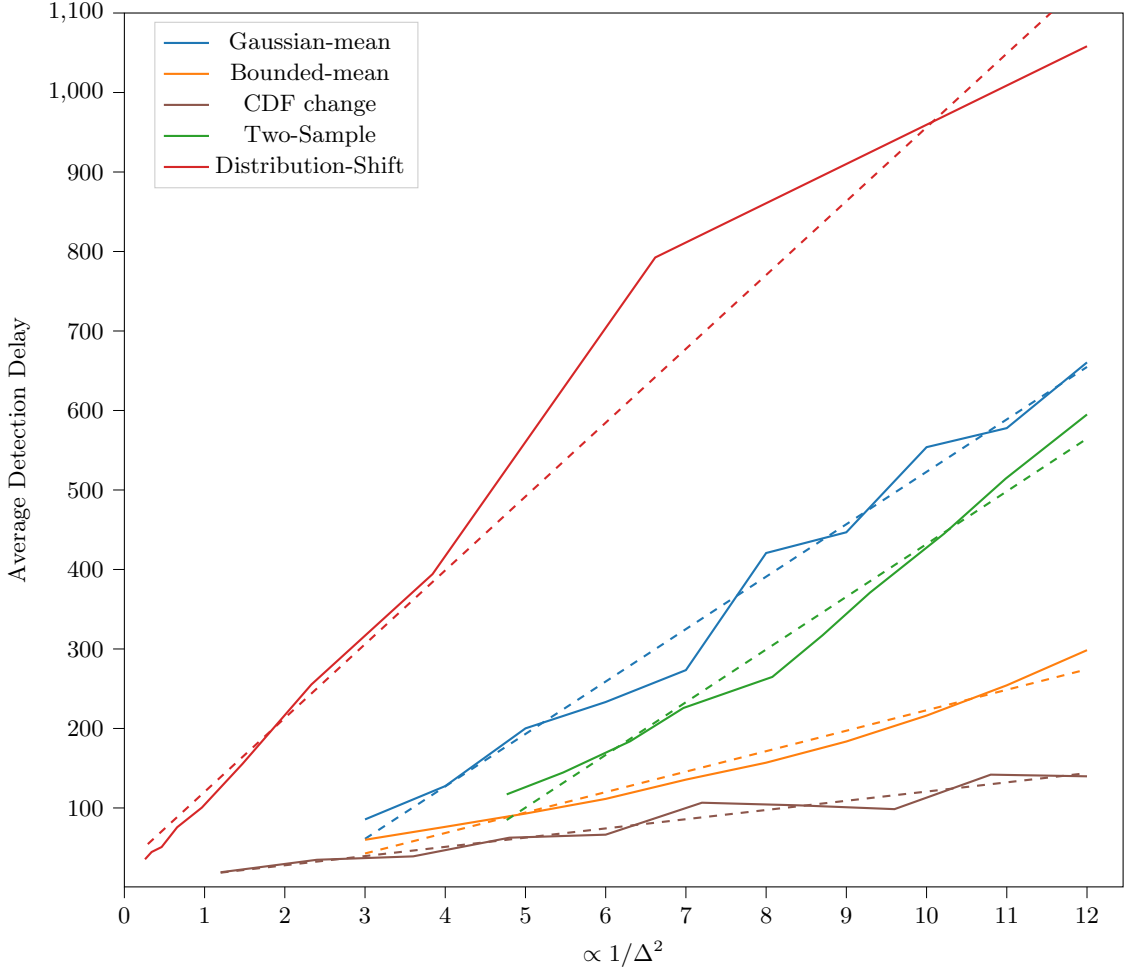
\begin{figure}[htb]
         \def\figwidth{0.9\linewidth}
             \def\figheight{0.8\linewidth} %
             \hspace{-1em}
            \centering
            \input{Figures/GaussianMean/Delay_vs_Delta8208}
            \caption{In all the  instantiations of the \bcsdetector within \Cref{sec:instantiations},  the width of the CS  at a $\sqrt{(\log \log (t) + \log(1/\alpha)/t}$ rate. Hence, by~\Cref{theorem:general-strategy}, we expect the detection delay to have an inverse quadratic dependence on the magnitude of change~($\Delta$). The figure above verifies this claim empirically. 
           }
            \label{fig:Delay-vs-Delta}
    \end{figure}
    \fi

\subsection{Nonparametric Change of Mean Detection}
\label{subsec:bounded-mean} 
    We now consider a nonparametric analog of the change of mean detection problem from the previous section. Here we assume that $X_1, X_2, \ldots $ are independent random variables taking values in a bounded interval $\mc{X} \subset \mbb{R}$, which we set to $[0,1]$ without loss of generality. Prior to the changepoint $T$, we assume that the observations have a mean $\theta_0 \in \Theta = [0,1]$, while it changes to $\theta_1 \neq \theta_0$ after time $T$.
    In this case, we can use the empirical Bernstein~(EB) confidence intervals developed by~\citet{waudby2023estimating}. To state the closed-form expression of the EB confidence sequence, we first need to introduce the following terms:
        $\muhat_t = \frac{\frac{1}{2} + \sum_{i=1}^t X_i}{t+1}$,  $\sigmahat_t^2 = \frac{\frac{1}{4} + \sum_{i=1}^t (X_i - \muhat_t)^2 }{t+1}$,   $\lambda_t = \sqrt{ \frac{2 \log(2/\alpha)}{\sigmahat_t^2 t \log(t+1)}} \wedge \frac{1}{2},$ 
       $\thetahat_t = \frac{\sum_{i=1}^t \lambda_i X_i}{\sum_{i=1}^t\lambda_i}$, $v_t = (4/\log(2/\alpha)) (X_t - \muhat_{t-1})^2$,  and $\Psi_E(x) = \frac{- \log(1-x) - x }{4}$, for all  $x\in [0, 1)$. 
    Using these terms, we can now state the EB-CS derived by~\citet{waudby2023estimating} as follows: 
        $C_t = [ \thetahat_t + w_t/2, \thetahat_t - w_t/2 ]$, 
    with  
        $w_t = \frac{ \log(2/\alpha) + \sum_{i=1}^t v_i \Psi_E(\lambda_i) }{\sum_{i=1}^t \lambda_i}$. 
    Again, by an application of the general result,~\Cref{theorem:general-strategy}, we can get the following bound on the expected detection delay of the SCD scheme, that uses the above EB-CS. 
    \begin{proposition}
        \label{prop:bounded-mean}
        Suppose $X_1, X_2, \ldots, X_{T}$ are drawn \iid from a distribution on $\mc{X} = [0,1]$ with mean $\theta_0 \in \Theta = [0,1]$, while $X_{T+1}, \ldots$ are drawn from $P_{\theta_1}$ with mean $\theta_1 \neq \theta_0$. Then, assuming $\theta_0$ is known, this instance of \bcsdetector~(\Cref{def:general-strategy}) satisfies the following~(with $\Delta \defined |\theta_0-\theta_1|$, and $\sigma_1^2 = \mathbb{E}_{P_{\theta_1}}[(X-\theta_1)^2]$): 
        \begin{align}
            \mathbb{E}_{T}[(\tau-T)^+] = \mc{O}\lp \sigma_1^2\, \frac{\log(1/\Delta) + \log(1/\alpha)}{\Delta^2} \rp. 
        \end{align}
    \end{proposition}
    
    \begin{remark}
        \label{remark:bounded-mean}
        While we stated the above result under the assumption that the pre- and post-change observations are \iid from $P_{\theta_0}$ and $P_{\theta_1}$ respectively, we note that similar results can be obtained when the random variables are only independent, with fixed (pre- and post-change) means. Further, the assumption that $\theta_0$ is known can also be waived, at the cost of conditioning on the `good' event $\mc{E}$. 
    \end{remark}

    \begin{remark}
        \label{remark:other-settings} In this subsection, we have considered the task of change-of-mean detection in perhaps the simplest (but nontrivial) nonparametric setting. The same ideas developed here, can however, we extended easily to other interesting cases, such as the sub-Gaussian family using the CS derived by~\citet{howard2021time}, or for heavy-tailed distributions~\citep{wang2022catoni}. 
    \end{remark}

    \noindent\textbf{Empirical Verification.} To verify our theoretical claims, we consider the case of an independent stream of observations supported on $\mc{X} = [0,1]$, with pre- and post-change parameters, $\theta_0$ and $\theta_1$ respectively. We define distributions with specified means by taking approprirate mixtures of uniform distributions, as described in~\Cref{appendix:experiments}.   We plot the performance of our changepoint detection scheme for a fixed problem instance with $(\theta_0, \theta_1, T)= (0.4, 0.6, 800)$ in~\Cref{fig:Bernoulli-mean-CS} (\Cref{appendix:experiments}). The predicted inverse quadratic dependence of the detection delay is verified in~\Cref{fig:Delay-vs-Delta}. 
    
\subsection{Detecting Changes in CDFs}    
\label{subsec:cdfs}
    Staying with real-valued observations~(or more generally, observations on totally ordered spaces), we now consider a more general question of detecting whether there have been any changes in distribution generating the observations. Since real valued random-variables are completely characterized by their cumulative distribution functions~(CDFs), this task can be framed in terms of detecting changes in the CDFs. More formally, we assume that we are given a stream of observations, $X_1, X_2, \ldots$, that are drawn according to a distribution $\theta_0 = F_0$ for $t < T$, and according to a distribution $\theta_1 = F_1$ for $t \geq T$. Thus, in this case, $\Theta$ is the infinite-dimensional space of all feasible CDFs on $\mbb{R}$, and we endow it with the Kolmogorov-Smirnov~(KS) metric, $\dks$, defined as $\dks(F, G) = \sup_{x \in \mbb{R}} \; |F(x)-G(x)|$. 
    
    To instantiate our SCD scheme, we will employ the following level-$\alpha$ confidence sequence for the CDF in terms of the  KS metric, recently derived by~\citet{howard2022sequential}: $C_t = \{\theta \in \Theta: \dks(\theta, \thetahat_t) \leq w_t/2\}$,  
    where $ w_t = 1.7 \sqrt{ \log \log(et)  + 0.8 \log(1612/\alpha)/ t}$. 
    As a consequence of the general result in~\Cref{theorem:general-strategy}, we can obtain the following performance guarantee for this scheme. 
    \begin{corollary}
        \label{prop:cdf}
        Suppose, for some $T < \infty$, the observations $X_1, \ldots, X_{T}$ are drawn from a known distribution $F_0$, while for $t \geq 1$, the observations $X_{T+1}, X_{T+2}, \ldots$ are drawn from an unknown $F_1$. Then, our SCD scheme instantiated with the CS stated above satisfies~(with $\Delta \defined \dks(F_1, F_0)$):
        \begin{align}
            \mathbb{E}_T[(\tau - T)^+] = \mc{O} \lp  \frac{ \log \log(1/\Delta) + \log(1/\alpha)}{\Delta^2}\rp. 
        \end{align}
    \end{corollary}

    \noindent\textbf{Empirical verification.} We test the performance of our proposed scheme for $t$-distributions with $3$ degrees of freedom. In~\Cref{fig:CDF-CS} in~\Cref{appendix:experiments}, we show the performance of our proposed scheme for a fixed problem instance where the pre- and post-change CDFs satisfy $\Delta = \dks(F_0, F_1) \approx 0.4$. The variation of the average detection delay with changing values of $\Delta$ is plotted in~\Cref{fig:Delay-vs-Delta}, and it displays the expected inverse quadratic dependence.

\subsection{Detecting Change in Homogeneity of Two Streams}
\label{subsec:two-sample}
    Suppose we have a stream of observations  in a product space $\mc{X} = \mc{U} \times \mc{U}$, and the parameter set consists of all product distributions on $\mc{U}\times \mc{U}$; that is, $\Theta = \{ P\times Q: P, Q \in \mc{P}(\mc{U})\}$. Prior to the changepoint, we assume that the observations $X_1, X_2, \ldots, X_{T}$, with $X_t = (U_t, V_t) \in \mc{U} \times \mc{U}$, are drawn from $\theta_0 = P_U \times P_V$; while the post-change observations $X_{T+1}, X_{T+2}, \ldots$ are assumed to be drawn from some other product distribution $Q_{U} \times Q_V$. Given some statistical distance measure, $\rho: \mc{P}(\mc{U}) \times \mc{P}(\mc{U}) \to \mbb{R}$, we assume that the $\rho(P_U, P_V) \neq \rho(Q_U, Q_V)$. An interesting special case of this problem, motivated by the two-sample  testing problem, is when $P_U = P_V$, and $Q_U \neq Q_V$.  
    
    If $\rho$ is a probability metric, we can use it induce a distance metric, $d$,  on the parameter space $\Theta$ as follows: $d(\theta_0, \theta_1) = |\rho(P_U, P_V) - \rho(Q_U, Q_V)$, where  $\theta_0 = P_U \times P_V$ and $\theta_1 = Q_U \times Q_V$.  Then, to obtain a changepoint detection scheme we can employ  CSs for the statistical distance $\rho$. We instantiate this strategy with the kernel-MMD metric~(defined in~\Cref{appendix:background}) associated with a kernel $k$, denoted by $\dmmd(\cdot, \cdot)$. 
    We use the following CS derived by \citet{manole2021sequential} for the kernel-MMD  distance between two distributions~(assuming that $\sup_{u, u'} k(u, u') \leq 1$):
        $C_t = \{(P, Q): \dmmd(\Phat_t, \Qhat_t) - \gamma_t \leq  \dmmd(P, Q) \leq     \dmmd(\Phat_t, \Qhat_t) + 2\kappa_t \}$ where
        $\kappa_t = \sqrt{ \big(\log ( (1 \vee \log_2(t))^2 \pi^2/6) + \log(4/\alpha)\big)/t}$, and  $\gamma_t = \big(4\sqrt{2}/\sqrt{t}\big) \lp 1 + \sqrt{ \log \lp 3.54 e(1 \vee \log_2t)^3 \rp + \log(2/\alpha)}\rp.$
   For this instantiation,~\Cref{theorem:general-strategy} implies  the following.
    \begin{corollary}
        \label{prop:kernel-MMD}
        With $X_t$ denoting the pair $(U_t, V_t)$ on $\mc{X} \times \mc{X}$, suppose that $X_1, \ldots, X_{T}$ are drawn \iid from $\theta_0 = P_U \times P_V$, and $X_{T+1}, X_{T+2}, \ldots$ are drawn \iid from a distribution $ \theta_1 = Q_U \times Q_V$. Then, with $\Delta > 0$ denoting  $d(\theta_0, \theta_1) = |\dmmd(P_U, P_V) - \dmmd(Q_U, Q_V)|$, we have the following upper bound on the expected detection delay of our BCS-detector based on the CS described above: 
        \begin{align}
            \mathbb{E}\lb (\tau-T)^+ | \mc{E}\rb = \mc{O} \lp \frac{ \log(1/\alpha) + \log \log(1/\Delta)}{\Delta^2} \rp, 
        \end{align}
        where $\mc{E}$ denotes the `good' event $\{ \dmmd(P_U, P_V) \in C_t:  t \leq T-1\}$ associated with the forward CS $\{C_t: t \geq 1\}$. 
    \end{corollary}
    
    \begin{remark}
        Consider the special case mentioned earlier, where $P_U=P_V=Q_U = P$ for some distribution $P$, and $Q_Y = Q \neq P$ for some other distribution. Furthermore, assume that it is known that prior to changepoint $P_U=P_V$~(that is, the event $\mc{E}$ is a probability one event). Then, the above result in this case implies an upper bound on the expected detection delay of $\mc{O}\lp \frac{ \log(1/\alpha) + \log \log(1/\Delta)}{\Delta^2} \rp$, with $\Delta = \dmmd(P, Q)$. This matches existing  results, such as the kernel CuSum scheme of \citet{wei2022online}. 
    \end{remark}
    
    \begin{remark}
        \label{remark:two-sample-distances} We focused on the case of the kernel-MMD metric mainly due to its generality (it is applicable to  distributions over arbitrary spaces on which positive definite kernels can be defined).  However, the same ideas are applicable to any statistical distance measure that is convex in its arguments, by using the reverse submartingale based confidence sequence construction of~\citet{manole2021sequential}. This family includes all  popular statistical distances, such as Wasserstein metrics, $f$-divergences and general integral probability metrics.
    \end{remark}
    \begin{remark}
        \label{remark:computational-cost} 
        The overall computational cost of our scheme is $\mathcal{O}(\tau^3)$, as our scheme involves constructing a new  backward CS, with $\mc{O}(t^2)$ cost, every round. 
        In practice, this complexity can be reduced, either by using linear or block-MMD statistics, and/or by  computing a new backward CS less frequently~(instead of doing so every round).  
    \end{remark}
    
    \noindent\textbf{Empirical Verification.} We study the performance of our scheme on  a stream of paired multivariate Gaussian observations in $p=5$ dimensions. The pre-change distributions, $P_U$ and $P_V$ both have zero mean and identity covariance; while for the post change distributions we have $Q_U = P_U$, and $Q_Y$ has a mean $\delta \boldsymbol{1}$, and a diagonal covariance matrix with randomly chosen values. 
     In~\Cref{fig:TwoSample-mean-CS} in~\Cref{appendix:experiments}, we plot the performance of our changepoint detection scheme for a fixed problem instance with $\Delta = \dmmd(Q_U, Q_V) \approx 0.33$, and $T=800$, while the inverse quadratic dependence of the average detection delay with $\Delta$ is verified in~\Cref{fig:Delay-vs-Delta}. 

\subsection{Detecting Harmful Distribution Shifts}
\label{subsec:harmful-shifts}
    As a final application, consider the task of detecting `harmful' changes between train and test distributions of a machine learning~(ML) model. Following~\citet{podkopaev2021tracking}, we are interested in detecting only those distribution changes that lead to a sufficiently large increase in the risk~(i.e., expected loss) of the trained ML model. 

    Formally, suppose a machine learning model, denoted by $h$,  is trained on a dataset drawn \iid from a \emph{source} distribution $P_S$ taking values on some space $\mc{X}$. For some bounded loss function, $\phi$, we let $\theta_0$ denote the expected training loss of this model; $ \theta_0 = \mathbb{E}_{P_S}[\phi(X, h)]$.
    Next, we assume that the model $h$ is deployed on a stream of test data, denoted by $X_1, X_2, \ldots$, drawn from the source (or training) distribution $P_S$ for $t<T$; and from some other distribution $P_T \neq P_S$ with $P_T \neq P_S$.  
    Our goal is to detect post-change distributions $P_T$ that are `harmful' to the trained model; that is, they result in  an increase in expected loss: $\theta_1 \defined \mathbb{E}_{P_T}[\phi(h, X)] > \theta_0$~(see~\Cref{fig:harmful-distribution-shift-1} in~\Cref{appendix:experiments}). 

    For bounded loss functions $\phi$,  this problem fits into the  nonparametric change of mean detection framework of~\Cref{subsec:bounded-mean}.  Since we are only interested in one-sided changes, we can modify the strategy of~\Cref{subsec:bounded-mean} to use only upper CS in the forward direction, and lower CSs in the backward direction. As in~\Cref{prop:bounded-mean}, for this strategy,  we can show that the  expected detection delay of the scheme will depend inversely on how `harmful' the target distribution is (i.e., the gap $\theta_1 - \theta_0$). 
    
    \noindent\textbf{Empirical Verification.} To illustrate the ideas discussed above, we consider a simple binary classification problem with linear classifiers and $2$-dimensional features~(see~\Cref{appendix:experiments} for details). We plot the performance of our scheme on a specific problem with $\Delta\approx 0.16$ in~\Cref{fig:Distribution-shift-CS} in~\Cref{appendix:experiments}, and also verify the inverse quadratic dependence of average detection delay on $\Delta$ in~\Cref{fig:Delay-vs-Delta}.

\subsection{Other change detection tasks}
\label{subsec:other-tasks} 
    We have illustrated the generality of our \bcsdetector strategy by instantiating it for five different scenarios in this section. For simplicity, we focused mainly on univariate observations~(with the exception of~\Cref{subsec:two-sample}). However, we note that the same ideas used in the previous instantiations also carry over easily to more general observations, or under additional robustness or privacy constraints. We list some such examples here, without going into the details of analysis or practical implementations: \\
        \textbf{(i) Exponential family.} In this case, the observations $X_1, X_2, \ldots$ lie in $\mc{X} = \mathbb{R}^p$, and the pre- and post-change distributions $(P_{\theta_0}$ and $P_{\theta_1}$) are chosen from a finite-dimensional exponential family with $\Theta = \mathbb{R}^{m}$ for some $m<\infty$. Here, we can use the \bcsdetector with the CSs derived by~\citet{chowdhury2022bregman}. \\
        \textbf{(ii) Covariance matrix.} Again, we assume that $\mc{X}=\mathbb{R}^p$, but now we assume that at the change point $T$, the covariance matrix of the observations changes from $\theta_0 \in \mathbb{R}^{p\times p}$ to some $\theta_1 \neq \theta_0$. For this problem, we can instantiate the \bcsdetector with the CS for covariance matrices derived by~\citet[\S~4.3]{howard2021time}. \\
        \textbf{(iii) Nonparametric regression.} Suppose $U_1, U_2, \ldots $ denote \iid uniform draws from $\mc{U}=[0,1]^p$, for some $p \geq 1$. Let $\Theta$ denote an RKHS (with kernel $k$) of functions from $\mc{U}$ to $\mathbb{R}$. Then, for any $\theta \in \Theta$, define the random variable $Y_t \equiv Y_t(\theta) = \theta(U_t) + \eta_t$, where $\{\eta_t: t \geq 1\}$ are an \iid sequence of $1$-sub-Gaussian noise. Clearly, the joint distribution of $(U_t, Y_t)$ is parametrized by $\theta$. Consider the SCD problem, where $\theta=\theta_0$ prior to changepoint $T$, and $\theta = \theta_1$ after that, with $\|\theta_0-\theta_1\|_k >0$. Our \bcsdetector strategy is easily applicable to this scenario, with an infinite-dimensional index set $\Theta$, using the CS constructed by~\citet{chowdhury2017kernelized}. \\
        \textbf{(iv) Robust SCD.} An interesting variant of the SCD problem involves detecting changepoints under adversarial contamination~\citep{li2021adversarially}. Our \bcsdetector strategy readily extends to such scenarios, by exploiting recent robust confidence sequence constructions, such as those by \citet{wang2023huber}. \\
        \textbf{(v) Private SCD.} Privacy is an important concern in many applications, especially involving personal data, and is often ensured by revealing only randomized versions of the actual data to the analyst. This adds another layer of complexity to the usual SCD task~\citep{cummings2018differentially}. However, our \bcsdetector framework can easily handle this, building upon the recent advances in private CS construction~\citep{waudby2022locally}. 
        
\section{Conclusion}
\label{sec:conclusion} 
    We proposed a general strategy~(\bcsdetector) for designing sequential changepoint detection~(SCD) schemes  by carefully combining confidence sequences~(CSs), and backward CSs --- a novel variant of CSs, that we introduced in this paper. Under very mild, and natural requirements on the CSs, we showed that \bcsdetector provides tight control over the ARL and the detection delay. Leveraging the recent progress in constructing CSs, we  instantiated our strategy for a wide range of SCD problems~(both parametric, and nonparametric), and empirically verified the theoretical claims via some small-scale numerical experiments. 
    
    Our work opens up several directions for future work:\\
        \textbf{(i)} Constructing a new backward CS in every round can be computationally costly, and in most cases results in an overall quadratic~(or even cubic) complexity. An interesting direction to pursue is to investigate if we can achieve the same performance by updating the backward CS fewer times. 
        \textbf{(ii)} In~\Cref{remark:changepoint-estimate} we defined estimators of the changepoint~($T$), and the change magnitude~($\Delta$), which performed well empirically as shown in~\Cref{fig:Gaussian-mean-CS} and~\Cref{appendix:experiments}. Establishing theoretical guarantees for them is another interesting question for future work. 
\newpage 
\bibliography{main}
\bibliographystyle{abbrvnat}

\newpage 
\onecolumn
\appendix
\section{Additional Background}
\label{appendix:background}

\begin{figure}[htb]
        \centering
        \begin{tabular}{ccc}
            \includegraphics[width=.48\textwidth]{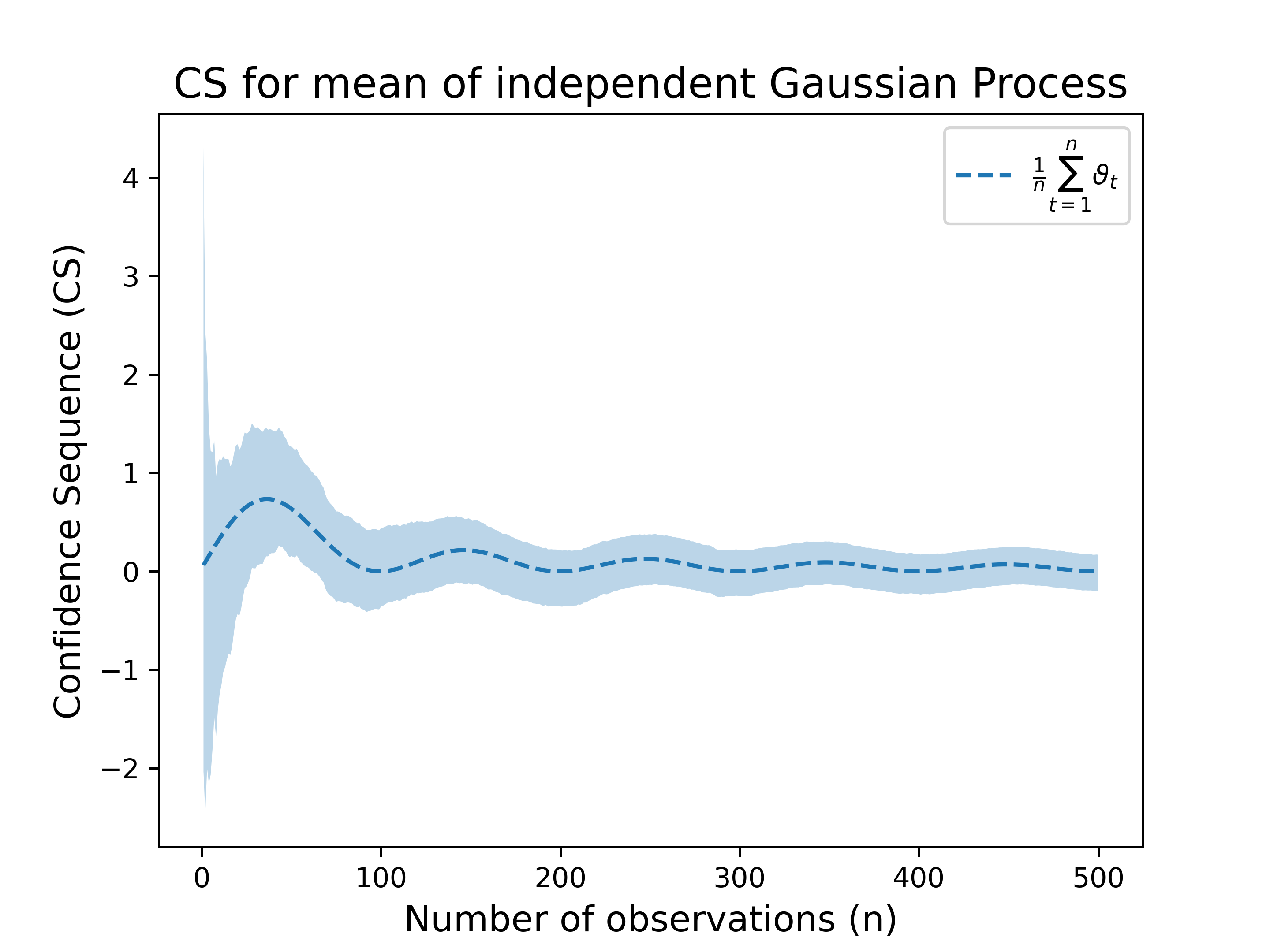}
        \end{tabular}
        \caption{An example of the CS introduced in~\Cref{example:gaussian-mean-cs} for the running (conditional) mean of independent Gaussian processes with a time-varying mean function, variance fixed at $1$.}
        \label{fig:gaussian-cs-example}
    \end{figure}

    \noindent\textbf{CS for non-\iid observations.} As mentioned in~\Cref{remark:CS-def-general}, we can also define CSs for non-\iid observations, as follows. 
    \begin{definition}
        \label{def:CS-def-general}
        Suppose $X_1, X_2, \ldots,$ denote an independent stream of observations on $\mc{X}$, with $X_t \sim P_{\vartheta_t}$ for $\vartheta_t \in \Theta$. 
        Assuming $\Theta$ is a vector space, we say that $\{C_t \subset \Theta: t \geq 1\}$ is a level $(1-\alpha)$ CS for the running average of parameters if 
        $\mathbb{P}\lp \{\forall t \geq 1: (1/t)\sum{i=1}^t \vartheta_t \in C_t \} \rp \geq 1-\alpha$.  
    \end{definition}
    We show a simple illustration of the CS introduced in~\Cref{example:gaussian-mean-cs} with time varying parameters in~\Cref{fig:gaussian-cs-example}. 
    
\begin{figure}
            \centering
            \includegraphics[width=0.3\linewidth]{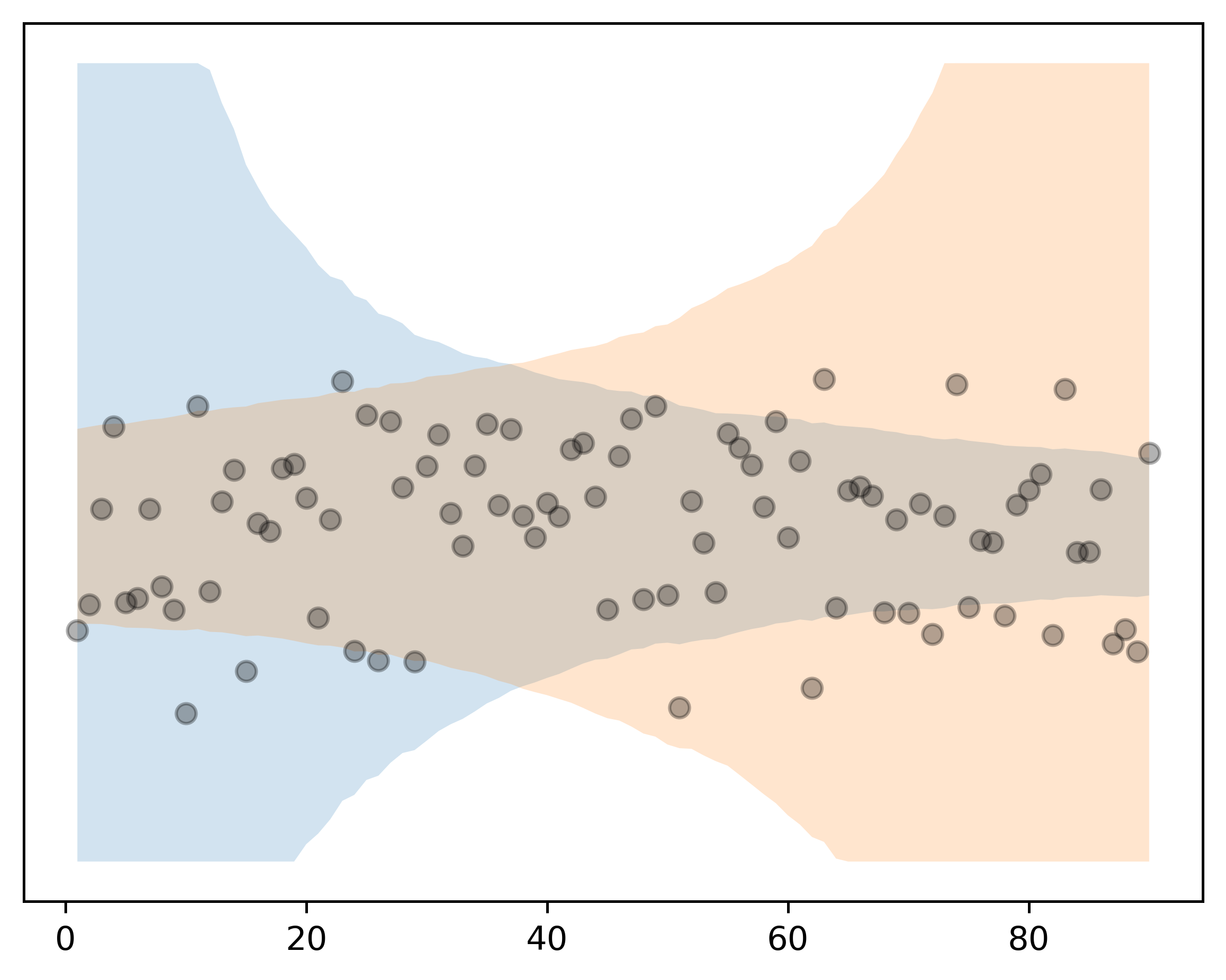}
            \includegraphics[width=0.3\linewidth]{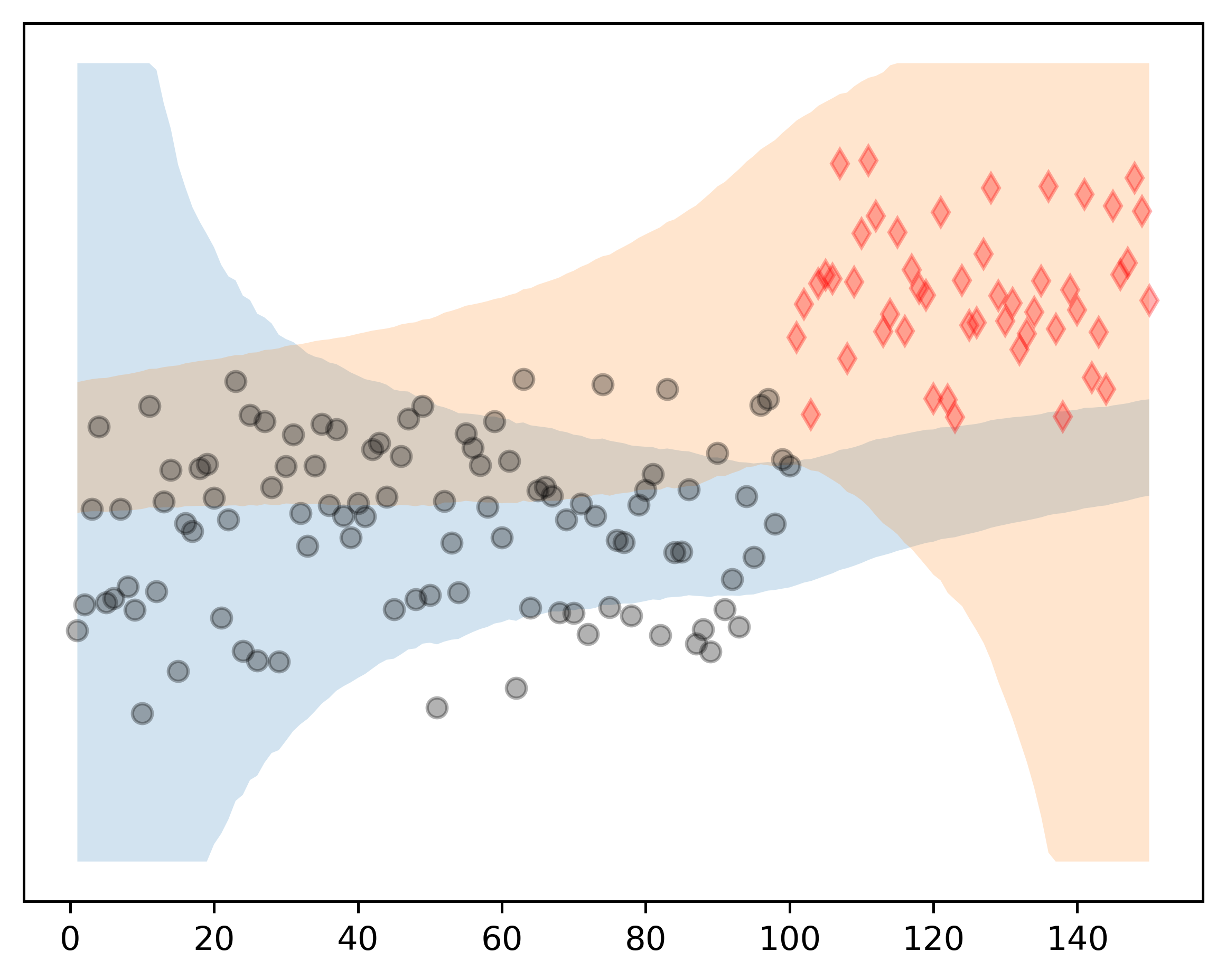}
            \includegraphics[width=0.3\linewidth]{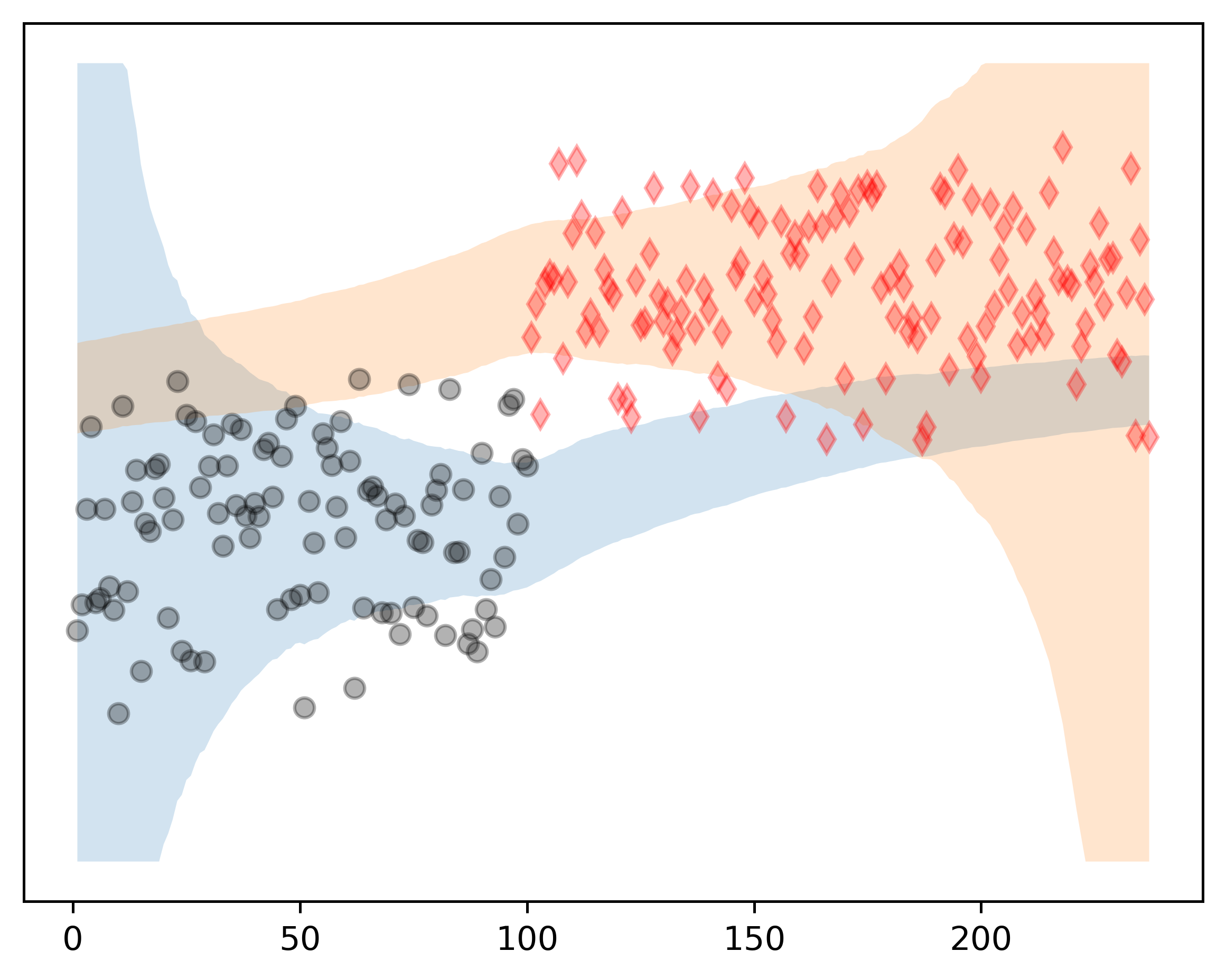}
            \caption{The plots illustrate the general ideas underlying our \bcsdetector strategy. Prior to the changepoint~(first plot), both the forward and backward CSs have significant overlap. After getting some post-change observations~(middle plot), the backward CS starts to drift away from the forward CS, although the deviation is not enough for the two CSs to become inconsistent. Finally, the last plot shows the scenario, where a sufficiently large number of post-change observations have arrived, which causes the forward and backward CSs to disagree. When this occurs, our scheme stops and rejects the null.}
            \label{fig:illustrating-bcsdetector}
        \end{figure}

    \noindent\textbf{Implementing backward CSs.} If we know how to construct forward CSs, we can use that directly to construct backward CSs in the following steps:  
    \begin{itemize}\itemsep0em
        \item  At the end of round $n$, we have observed $X_1, \ldots, X_n$. Introduce the time-reversed version of the observations $Y_s = X_{n+1-s}$ for $s \in [n] \defined \{1, \ldots, n\}$. \\
        \item  Construct a new level-$(1-\alpha)$ CS using $\{Y_s: s \in [n]\}$, denoted by $\{\barB_s^{(n)}: s \in [n]\}$. Note that $\barB_s^{(n)}$ is $\sigma(Y_1, \ldots, Y_s) = \sigma(X_n, X_{n-1}, \ldots, X_{n-s+1})$ measurable for all $s \in [n]$. \\
        \item  Finally, we again reverse the index of the CS, to obtain $\{B_t^{(n)}: t \in [n]\}$, where $B_t^{(n)}= \barB_{n+1-s}^{(n)}$. Note that by virtue of being a CS, we have  $\prob_{\infty} \big( \forall t \in [n]: \theta_0 \in B_t^{(n)} \big) \geq 1-\alpha$. 
        The superscript $(n)$ serves as a reminder that  there is only one forward CS, but there is a different backward CS constructed afresh at each time $n$.
    \end{itemize}
    \noindent\textbf{Illustration of \bcsdetector strategy.}
        In~\Cref{fig:illustrating-bcsdetector},  we illustrate the intuition underlying our general \bcsdetector strategy, introduced in~\Cref{sec:proposed-scheme-2}, using the task of detecting change in means of bounded observations~(\Cref{subsec:bounded-mean}).  The three plots in~\Cref{fig:illustrating-bcsdetector} highlight the following aspects of our scheme: 
        \begin{itemize}
            \item Prior to the changepoint~(or if there is no changepoint), both the forward and backward CSs concentrate around the same region in parameter space $\Theta$~(in this case, $[0,1]$). In particular, note that in this case, for all values of $t$, one of the confidence intervals~(CI) is entirely contained in the corresponding CI of the other CS. 
            \item As observations from the post change distribution start arriving, we expect the forward and backward CSs to start drifting away from each other. This is illustrated in the second plot of~\Cref{fig:illustrating-bcsdetector}. 
            \item Finally, after sufficiently many post-change observations, the backward CS starts concentrating around the post-change parameter $\theta_1$, as a result of which some of the backward CIs become disjoint with the forward CIs. Our \bcsdetector scheme uses this occurrence as a signal to stop, and declare a changepoint. 
        \end{itemize}
    \noindent\textbf{The kernel-MMD metric.} In~\Cref{subsec:two-sample}, we constructed a scheme for detecting changes the pairwise kernel-MMD distance between the distributions generating a stream of paired observations. Here, we recall the its definition. 
    
    We assume that $k:\mc{X} \times \mc{X} \to \mbb{R}$ denotes a uniformly-bounded positive-definite  kernel, and let $\mc{H}_k$ denote the  reproducing kernel  Hilbert space~(RKHS) associated with $k$.  
    \begin{definition}
        \label{def:kernel-MMD}
        Given a positive definite kernel $k:\mc{X} \times\mc{X} \to \mbb{R}$, the kernel-MMD distance between two-distributions $P$ and $Q$ on $\mc{X}$ is defined as 
        \begin{align}
            \dmmd(P, Q) = \sup_{g \in \mc{H}_k: \|g\|_k \leq 1} \; \expec_P[g(X)] - \expec_Q[g(Y)]. 
        \end{align}
    \end{definition}
    The kernel-MMD distance defined above is an instance of a class of statistical distances called \emph{integral probability metrics~(IPMs)}. For a class of kernels, called characteristic kernels, it is known that $\dmmd$ is a distance metric on the space of probability distributions. 

\section{Proof of~\Cref{prop:general-strategy-1}}
\label{appendix:proof-general-strategy-1}
    \textbf{Proof of the bound on probability of false alarm.} Recall that the stopping time is defined as the first time, $\tau$, at which we have the condition $\cap_{t=1}^{\tau}C_t = \emptyset$. Now, consider the `good' event of the CS under the null: $\mc{E} = \cap_{t=1}^{\infty}\{\theta_0 \in C_t\}$, which satisfies $\prob_\infty(\mc{E}) \geq 1-\alpha$ by definition. Hence, under this event $\{\theta_0\} \subset \cap_{t=1}^{\infty} \neq \emptyset$, which in turn implies that $\prob_\infty(\tau=\infty) \geq 1- \alpha$, as required. 
    
    \textbf{Proof of the bound on detection delay.} Let $w(t) \equiv w(t, \Theta, \alpha)$ denote the width of the confidence $C_t$ after $t$ observations. By assumption, $T$ is large enough to ensure that under the `good' event $\mc{E} = \cap_{t=1}^{\infty} \left \{ \frac{1}{t} \vartheta_t  \in C_t \right\}$, we have that $\theta_1 \not \in C_T$ at the changepoint. Note that in the definition of $\mc{E}$, we have $\vartheta_t = \theta_0\bone_{t \leq T} + \theta_1 \bone_{t > T}$. 
    
    Under the event $\mc{E}$, we know that $\theta_0  \in C_{T}$. For any $t> T$, introduce the terms $\lambda_t = T/t$ and $\bar{\lambda}_t = 1-\lambda_t$. Then, by definition of confidence sequences, we have $\lambda_t \theta_0 + \bar{\lambda}_t \theta_1 \in C_t$ for all $t>T$ under the event $\mc{E}$. The width of the set $C_t$ at $t>T$ is no larger than $w(t) \equiv w(t, \Theta, \alpha)$. Hence, a sufficient condition for stopping prior to $t>T$ is if the sum of the widths of $C_t$ and $C_T$, that is $w(t) + w(T)$, is smaller than $d\lp \theta_0, \lambda_t\theta_0 + \bar{\lambda}_t \theta_1\rp$. 
    
     \textbf{Informal calculations for~\Cref{remark:norm-induced}.} As mentioned in~\Cref{remark:norm-induced}, in many cases, the stopping time $\tau$ satisfies: $ \tau \approx \min \{t \geq T: 1/\sqrt{t} + 1/\sqrt{T} \leq ((t-T)/t) \Delta\}$, where $\Delta = d(\theta_0, \theta_1)$. We now consider the behavior of the delay, $\tau - T$, in two different regimes of the changepoint $T$. 
     
     First we consider the case where $T$ is `small'; that is $T \approx 1/\Delta^2$. For concreteness, assume that $T = 9/\Delta^2$. Then, it is easy to check that with $\tau \leq 4T$, since for $t = 4\tau$
     \begin{align}
         \frac{1}{\sqrt{t}} + \frac{1}{\sqrt{T}} = \frac{\Delta}{4} + \frac{\Delta}{8} \leq \frac{3\Delta}{4} = \frac{t - T}{t} \Delta. 
     \end{align}
     Next, we consider the case where $\Delta$ is fixed, but $T \to \infty$. In this case, we have $\tau - T = \mc{O}\lp \sqrt{T} \rp$. To see this, consider $t = T + u$, with $u = o(T)$. Then, we have 
     \begin{align}
         \frac{1}{\sqrt{t}} + \frac{1}{\sqrt{T}} \approx \frac{2}{\sqrt{T}}, \quad \text{and} \quad \frac{t-T}{t}\Delta \approx \frac{u}{T} \Delta. 
     \end{align}
     Thus, this implies that the appropriate order of growth of the detection delay is $u = \mc{O}(\sqrt{T}/\Delta)$, as $T \to \infty$ with $\Delta$ fixed. 
     
     \begin{figure}[htb!]
         \centering
         \includegraphics{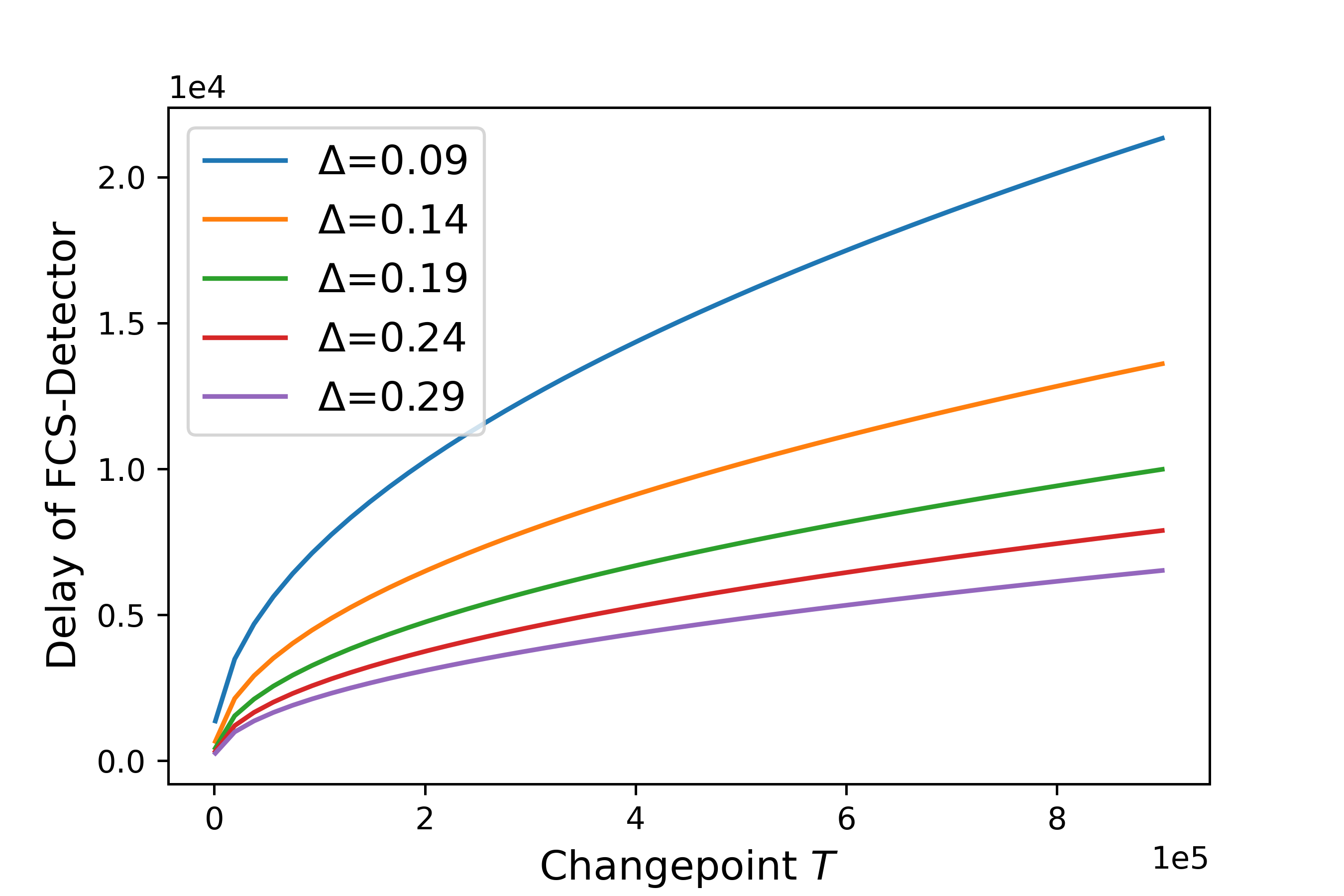}
         \caption{The plots show the variation of the delay of SCD scheme introduced in~\Cref{sec:proposed-scheme-1}, under the conditions of~\Cref{remark:norm-induced}. When the changepoint $T \approx 1/\Delta^2$, the delay has a linear dependence on $T$, while in the regime where $T \to \infty$ with $\Delta$ fixed, the delay behaves  $\approx \sqrt{T}$.}
         \label{fig:delay-strategy-1}
     \end{figure}

\section{Proof of~\Cref{theorem:general-strategy}}
\label{appendix:proof-general-strategy}

        \textbf{Proof of the ARL control.}
       To prove this result, note that for us to stop under the null at some time $\tau$, either the forward or the backward CS (or both) must be miscovering and failing to contain $\theta_0$. 
        In other words,
        \begin{align}
        \{\tau = N, T=\infty\}  &\implies \{C_N \text{ miscovers}\} \vee  \{ B^{(N)}_1 \text{ miscovers} \}\\
        \{\tau \leq N, T=\infty\}  &\implies \{(C_t)_{t=1}^N \text{ miscovers}\} \vee \bigcup_{n=1}^N \{(B^{(n)}_s)_{s=n}^1 \text{ miscovers} \}
        \end{align}
        Since each is a $(1-\alpha)$-CS, we have by a simple union bound that for any fixed time $N$,
        \[
        \text{Pr}_\infty(\tau \leq N) \leq \min\{(N+1)\alpha, 1\}.
        \]
        Rephrasing, we have $\Pr_\infty(\tau > N) \geq (1-\alpha(N+1)) \vee 0$, and thus
        \begin{align}
        \EE_\infty \tau &= \sum_{N=1}^\infty \text{Pr}_\infty(\tau > N) \geq \sum_{N=1}^{1/\alpha -1} (1-\alpha-\alpha N)\\
        &= (1/\alpha - 1)(1-\alpha) - \alpha \sum_{N=1}^{1/\alpha -1} N\\
        &= 1/\alpha - 2 + \alpha - \alpha \frac{(1/\alpha-1)(1/\alpha)}{2} \\
        &= \frac1{2\alpha} -3/2 + \alpha.
        \end{align}
        
        As we alluded to earlier in~\Cref{remark:same-class}, we did not need to know $\theta_0$ or $\theta_1$ or the `direction' of the changepoint. The aforementioned argument goes through for any indexed family of distributions. 
        \bigskip
        
        \noindent \textbf{Proof of the detection delay bound.} We next obtain the upper bound on the average detection delay conditioned on the event $\mc{E} \defined \{ \theta_0 \in C_t: 1 \leq t \leq T\}$, stated in~\Cref{theorem:general-strategy}. Since, $(\tau - T)^+ = \max \{0, \tau - T\}$ is a non-negative random variable, we have 
        \begin{align}
            \mathbb{E}_{T}[(\tau - T)^+|\mc{E}] &= \sum_{t=0}^{\infty} \prob_T\lp (\tau-T)^+ \geq t |\mc{E} \rp \\
            & \leq \sum_{t=0}^{t'-1} 1  + \sum_{t \geq t'} \prob_T \lp (\tau-T)^+ \geq t' |\mc{E}\rp  \\
            & = \sum_{t=0}^{t'-1} 1  + \sum_{t \geq t'} \prob_T \lp \tau > T+t' |\mc{E}\rp  \\
            &= t' + \sum_{t \geq t'} \prob_T \lp \tau > T+t' | \mc{E}\rp, \quad \text{for any } t' \geq 1. \label{eq:delay-upper-bound-0}
        \end{align}
        Recall that $u_0 \equiv u_0(\theta_0, \theta_1, T)$ was introduced in the statement of~\Cref{theorem:general-strategy}, and it represents an upper bound on the smallest time after $T$ at which the  backward CS must stop intersecting with the forward CS. 
        For some integer $i\geq 1$, consider the event $\{\tau > T + i u_0\}$, and note that it satisfies the following inclusion: 
        \begin{align}
            \{\tau > T + i u_0\} \cap \mc{E} &\subset \lp  \cap_{j=1}^{i} \left\{ \theta_0  \in B_T^{(T+ju_0)} \right\} \rp \cap \mc{E}\\
            & \subset \lp  \cap_{j=1}^{i} \left\{ \theta_0  \in B_{T+(j-1)u_0}^{(T+ju_0)} \right\} \rp \cap \mc{E}. \label{eq:delay-upper-bound-2}
        \end{align}
        In the display above, \eqref{eq:delay-upper-bound-2} uses the fact that  $B_{s}^{(n)} \subset B_{s'}^{(n)}$ for any $s'>s$. Next, we make the following two observations: 
        \begin{itemize}
            \item [\textbf{(I)}] For $j \neq j'$, the events $\{\theta_0 \in B_{T+(j-1)u_0}^{(T+ju_0)}\}$ and  $\{\theta_0 \in B_{T+(j'-1)u_0}^{(T+j'u_0)}\}$ are independent 
            \item [\textbf{(II)}] For all $1\leq j \leq i$, the event $\{\theta_0 \in B_{T+(j-1)u_0}^{(T+ju_0)}\}$ is independent of $\mc{E}$. 
        \end{itemize}
       \sloppy The statement \textbf{(I)} follows from the observation that the event $\{\theta_0 \in B_{T+(j-1)u_0}^{(T+ju_0)}\}$  lies in the sigma-algebra $\sigma\lp \{X_k: T+(j-1)u_0 \leq k <T+ ju_0 \} \rp$, while  $\{\theta_0 \in B_{T+(j'-1)u_0}^{(T+j'u_0)}\}$ lies in $\sigma\lp \{X_k: T+(j'-1)u_0 \leq k <T+ j'u_0 \} \rp$; which are independent. Similarly, the second statement \textbf{(II)} uses the fact that $\mc{E}$ lies in $\sigma \lp \{X_k: 1 \leq k <T\} \rp$, which is independent of  $\sigma\lp \{X_k: T+(j-1)u_0 \leq k <T+ ju_0 \} \rp$ for all $1 \leq j \leq i$. 
       
        Based on the above observations, we conclude that 
        \begin{align}
            \prob_T \lp \tau > T + iu_0 | \mc{E} \rp & = \prob_T \lp \{ \tau > T+iu_0\} \cap \mc{E} \rp /  \prob_T(\mc{E})  \\
            &\leq  \prob_T\lp \lp \cap_{j=1}^i \{\theta_0 \in B_{T+(j-1)u_0}^{(T+ju_0)} \} \rp \cap   \mc{E} \rp / \prob_T(\mc{E})\\
            & =  \prob_T\lp \cap_{j=1}^i \{\theta_0 \in B_{T+(j-1)u_0}^{(T+ju_0)} \} \rp \label{eq:delay-upper-bound-5}\\
            & \leq  \alpha^{i}, \quad \text{for all } i \geq 1. \label{eq:delay-upper-bound-3}
        \end{align}
        In the above display,~\eqref{eq:delay-upper-bound-5}  follows from \textbf{(II)}, and~\eqref{eq:delay-upper-bound-3} uses~\textbf{(I)}. 
        Now, we return to~\eqref{eq:delay-upper-bound-0}, and set $t'$ to $i_0 \times u_0$ for some integer $i_0$ to be specified later. We then note that  
        \begin{align}
            \mathbb{E}_T[(\tau-T)^+ | \mc{E}] &\leq i_0 u_0 + \sum_{i=i_0}^{\infty} \sum_{t= i u_0}^{(i+1)u_0-1} \prob_T(\tau > T + t |\mc{E}) \\
            & \leq i_0 u_0 + \sum_{i=i_0}^{\infty} u_0 \prob_T \lp \tau > T + iu_0 |\mc{E} \rp \\ 
            & \leq i_0 u_0 + u_0 \frac{ \alpha^{i_0} }{1 - \alpha} = u_0 \lp i_0 + \frac{\alpha^{i_0}}{1-\alpha} \rp \label{eq:delay-upper-bound-4}. 
        \end{align}
        The inequality in~\eqref{eq:delay-upper-bound-4} uses the fact that $\prob_T(\tau > T + iu_0 ) \leq \alpha^{i}$ derived in~\eqref{eq:delay-upper-bound-3}. The final result, as stated in~\Cref{theorem:general-strategy}, then follows by selecting $i_0 = \lceil \log(1/1-\alpha) / \log(1/\alpha) \rceil$. Note that when $\alpha<0.5$, we have $i_0=1$. 
    
\section{Deferred proofs from~\Cref{sec:instantiations}}
\label{appendix:proof-of-instatiations}
    
    \subsection{Proofs of~\Cref{prop:gaussian-mean},~\Cref{prop:cdf},~\Cref{prop:kernel-MMD}}
        All these three results can be obtained as a direct consequence of the following proposition. 
        \begin{proposition}
            \label{prop:stopping-time-calculation}
            For some $\Delta>0$, define the time $t_0$ as 
            \begin{align}
                t_0 = \min \left \{ t \geq 1: c \sqrt{ \frac{ \log \log t + \log(1/\alpha)}{t}} \leq \frac{\Delta}{2} \right\}, 
            \end{align}
            where $c > 0$ is some constant. Then, we have 
            \begin{align}
                t_0= \mc{O}\lp c^2\, \frac{ \log \log(c/\Delta) + \log(1/\alpha)}{\Delta^2} \rp. 
            \end{align}
        \end{proposition}
        \begin{proof}
            Without loss of generality, we assume that $c=1$; or equivalently, we can replace $\Delta$ with $\Delta/c$. Now, note that  we can upper bound $t_0 \leq t_1 + t_2$, where 
            \begin{align}
                t_1 = \min \{ t \geq 1: \sqrt{ \log \log (t) / t }  \leq \Delta/4\}, \quad \text{and} \quad t_2 = \min \{ t \geq 1: \sqrt{\log(1/\alpha)/t} \leq \Delta/4\}. 
            \end{align}
            By a simple calculation, we can obtain $t_2 = \mc{O} \lp \log(1/\alpha)/ \Delta^2 \rp$. Hence to complete the proof, we will show that $t_1 = \mc{O}\lp \log \log (1/\Delta)/ \Delta^2 \rp$.  We proceed in two steps: (i) first we show that $t_1 \leq 32/\Delta^3$, and (ii) using this, we refine the result to show that $t_1 = \mc{O} \lp \log \log(1/\Delta) / \Delta^2 \rp$. 
            
            Let $t_3 = 32/\Delta^3$ for $\Delta\leq 1$. Then, observe that 
            \begin{align}
                \lp \frac{4}{\Delta} \sqrt{ \frac{ \log \log (t_3)}{t_3} } \rp^2 = \frac{16}{\Delta^2} \times \frac{ \Delta^3 \, \log \log (32/\Delta^3)}{32} = \frac{1}{2} \frac{ \log \log (32/\Delta^2)}{ 1/\Delta} \leq 0.63 < 1.  
            \end{align}
            The last inequality, along with the definition of $t_1$ implies that $t_1 \leq t_3$.
            
            Hence, $\log \log (t_1) \leq \log \log (t_3) = \log \log (32/\Delta^3) \leq 2.35 +  \log \log(1/\Delta)$. Thus, we have 
            \begin{align}
                \frac{ \log \log (t_1) }{t_1} \leq \frac{2.35 +  \log\log(1/\Delta)}{t_1}, 
            \end{align}
            which implies that $t_1 \leq \frac{16}{\Delta^2} \lp  2.35 + \log \log (1/\Delta) \rp = \mc{O} \lp \log \log(1/\Delta)/ \Delta^2 \rp$. 
            Combining this bound on $t_1$, with the previously obtained upper bound on $t_2$; and using the fact that $t_0 \leq t_1 + t_2$,  we get the required result. 
        \end{proof}
    
    \subsection{Proof of~\Cref{prop:bounded-mean}}
    \label{appendix:proof-bounded-mean}
        To prove the variance adaptive bound on the expected detection delay, we first show the following result for the width of the CS derived by~\citet{waudby2023estimating}. 
        \begin{proposition}
            \label{prop:width-bounded-CS} We can modify the backward CS used instantiating the \bcsdetector in~\Cref{subsec:bounded-mean}, to obtian a level $(1-2\alpha)$-backward CS (for every $n \geq T$), denoted by $\{B_t^{(n)}: 1 \leq t \leq n\}$, such that the width of $B_t^{(n)}$ satisfies 
            \begin{align}
                \label{eq:waudby-smith-width}
                w(t, \theta_1, \alpha) = \mc{O} \lp \frac{\sigma_1}{\sqrt{t}} \lp \sqrt{\log(1/\alpha)} + \sqrt{\log t} \rp \rp, \quad \text{for all } T \leq t \leq n. 
            \end{align}
        \end{proposition}
        \begin{remark}
            The main benefit of this result is that it characterizes the width of the backward CS (for $n \geq T$) explicitly in terms of the standard deviation~($\sigma_1$) of the post-change distribution $P_{\theta_1}$, unlike the original CS described in~\Cref{subsec:bounded-mean}, whose width depends on empirical estimates of $\sigma_1$. As a consequence of this result, we obtain~\Cref{prop:bounded-mean} by first appealing to~\Cref{corollary:known-pre-change}, and then repeating the calculations used to obtain~\eqref{prop:stopping-time-calculation}. 
        \end{remark}
        \begin{proof}
            Since we are only interested in characterizing the order with which the width of the CS decays~(and not the exact constants), we will not track the constants in our argument for this proof. In particular, we will use $A \lesssim B$ to indicate that by $A/B = \mc{O}(1)$, and  $A \approx B$ to indicate that $A \lesssim B$ and $B \lesssim A$. 
            
            We proceed in the following steps: 
            \begin{itemize}
                \item First, we show that we can construct a level-$(1-\alpha)$ confidence sequence for the empirical variance based on  samples from the post-change distribution. In particular, let $\sigmahat_n^2 = (\frac{1}{4} + \sum_{t=1}^n (X_t - \muhat_t)^2)/(n+1)$, with $X_1, X_2, \ldots \sim P_{\theta_1}$ \iid. Then, we have the following: 
                \begin{align}
                    \prob \lp \mc{E}_1 \rp \geq 1-\alpha, \quad \text{where} \quad \mc{E}_1 \defined \cap_{n \geq 1}  \left\{ |\sigmahat_n^2 - \sigma_1^2| = \mc{O}\lp \sqrt{\frac{\log \log n + \log (1/\alpha)}{n}} \rp \right \}.   \label{eq:variance-cs}
                \end{align}
                Thus, using this event, for any $n \geq T$, we can modify the backward CS to get a level-$(1-2\alpha)$ Backward CS, in which $\sigmahat$ is replaced by $\sigma_1$~(plus a small approximation error term) for $T \leq t \leq n$. In the next two steps, we characterize the width of these level-$(1-2\alpha)$ CSs.
                
                \item  Next, we  show that under the event $\mc{E}_1$, we have 
                \begin{align}
                    \sum_{i=1}^n v_i \psi_E(\lambda_i) \approx  \log(\sigma_1^2 n). \label{eq:psi-approx}
                \end{align}
                \item Under the same event, we then show that 
                \begin{align}
                    \frac{1}{\sqrt{n}} \sum_{t=1}^n \lambda_i \approx \sqrt{ \frac{ 2 \log (2/\alpha)}{\sigma_1^2}}. 
                    \label{eq:lambda-approx}
                \end{align}
            \end{itemize}
            Combining these results, we get that the width of the CS is of the order 
            \begin{align}
                w_n \approx \frac{1}{\sqrt{n}}  \frac{ \log(2/\alpha) + \sum_{t=1}^n v_t \psi_E(t)}{\frac{1}{\sqrt{n}} \sum_{t=1}^n \lambda_t } \approx \frac{\sigma_1}{\sqrt{n}} \lp \sqrt{\log(2/\alpha)} + \sqrt{\log(n)} \rp,  
            \end{align}
            as required. Thus, it remains to prove~\eqref{eq:variance-cs},~\eqref{eq:psi-approx}, and~\eqref{eq:lambda-approx}. 
            
            \textbf{Proof of~\eqref{eq:variance-cs}.} To prove this, we first introduce the usual unbiased estimate of the variance: $\sigmatilde_n^2 = \frac{1}{n(n-1)} \sum_{i=1}^n \sum_{j \neq i} \frac{(X_i - X_j)^2}{2}$.  Since the observations $X_1, X_2, \ldots \sim P_{\theta_1}$ are bounded, and lie in $[0,1]$, it is easy to verify that $|\sigmahat_n^2 - \sigmatilde_n^2| = \mc{O}(1/n)$, and hence $\sigmahat_n^2 \approx \sigmatilde_n^2$. 
            
            Since, $\sigmatilde_n^2$ is an instance of a U-statistic, and hence the process $\{\sigmatilde^2: n \geq 1\}$ is a reverse-martingale, adapted to the exchangeable filtration. Using this fact, along with the boundedness~(and hence sub-Gaussianity) of the random variable $\sigmatilde^2$ for all $n \geq 1$, we can use~\citet[Corollary 8]{manole2021sequential}, to conclude the time-uniform concentration result: $\mathbb{P}\lp \forall n \geq 1: |\sigmatilde^2 - \sigma_1^2| = \mc{O}(r_n) \rp \geq 1-\alpha$, where $r_n = \sqrt{\lp \log \log n + \log(1/\alpha) \rp/n}$.

            \textbf{Proof of~\eqref{eq:psi-approx}.} To show this, we recall the fact that $\psi_E(\lambda)/(\lambda^2/8) \to 1$ as $\lambda \to 0$. Hence, we have the following: 
            \begin{align}
                \sum_{i=1}^n v_i \psi_E(\lambda_i) &= \frac{4}{\log(2/\alpha)} \sum_{i=1}^n (X_i - \muhat_i)^2 \psi_E(\lambda_i) \approx  \frac{4}{\log(2/\alpha)} \sum_{i=1}^n (X_i - \muhat_i)^2 \frac{\lambda_i^2}{8} \\
                &\lesssim \frac{1}{\log(2/\alpha)} \sum_{i=1}^n (X_i - \muhat_i)^2 \frac{ \log (2/\alpha)}{i \sigmahat_{i-1}^2} \lesssim \frac{1}{\log(2/\alpha)} \sum_{i=1}^n (X_i - \muhat_i)^2 \frac{ \log (2/\alpha)}{i \sigmahat_{i-1}^2} \\
                &\approx \sum_{i=1}^n  \frac{(X_i-\muhat_{i-1})^2}{\sum_{j=1}^i(X_j-\muhat_{j-1})^2} \\
                & \leq \log \lp \sum_{i=1}^n (X_i - \muhat_{i-1})^2 \rp \label{eq:proof-psi-approx-1} \\
                & \approx \log \lp n \sigmatilde^2_n \rp  \approx \log \lp n \sigma_1^2 \rp. 
            \end{align}
            In the above display,~\eqref{eq:proof-psi-approx-1} follows by an application of the following lemma with $f(x) = 1/x$. 
            
            \begin{lemma}[\citet{orabona2019modern}, Lemma~4.13]
            Let $a_i \geq 0$ for all $i$, and $f:[0,\infty) \to [0,\infty)$ be an increasing function. Then 
            \[
            \sum_{t=1}^T a_t f(a_0+\sum_{i=1}^t a_i) \leq \int_{a_0}^{\sum_{t=0}^T a_t} f(x) dx.
            \]
            \end{lemma}
            This concludes the proof.
            
            \textbf{Proof of~\eqref{eq:lambda-approx}.} We proceed as follows with $r_i = \sqrt{ \big( \log \log i + \log(1/\alpha) \big)/i }$
            \begin{align}
                \frac{1}{\sqrt{n}} \sum_{i=1}^n \lambda_i &= \frac{1}{\sqrt{n}} \sum_{i=1}^n \sqrt{ \frac{ 2\log(2/\alpha)}{\sigmahat_i^2 \, i}}   \gtrsim \sqrt{ \frac{ 2 \log(2/\alpha)}{n}} \sum_{i=1}^n \frac{1}{\sqrt{i \sigma_1^2 (1 + r_i)}} \\
                & \gtrsim \sqrt{ \frac{ 2 \log(2/\alpha)}{ \sigma_1^2\,n}} \lp \sum_{i=1}^n \frac{1}{\sqrt{i }} -  \sum_{i=1}^n \frac{r_i}{\sqrt i} \rp \\
                & \approx \sqrt{ \frac{ 2 \log(2/\alpha)}{ \sigma_1^2\,n}} \times \sqrt{n} = \sqrt{ \frac{ 2 \log(2/\alpha)}{ \sigma_1^2}}. 
            \end{align}
            This completes the proof of~\Cref{prop:width-bounded-CS}. 
        \end{proof}

\section{Details of Experiments}
\label{appendix:experiments}

  \iftrue 
        \begin{figure}
             \def\figwidth{0.25\linewidth}
             \def\figheight{0.25\linewidth} %
            \centering
            \input{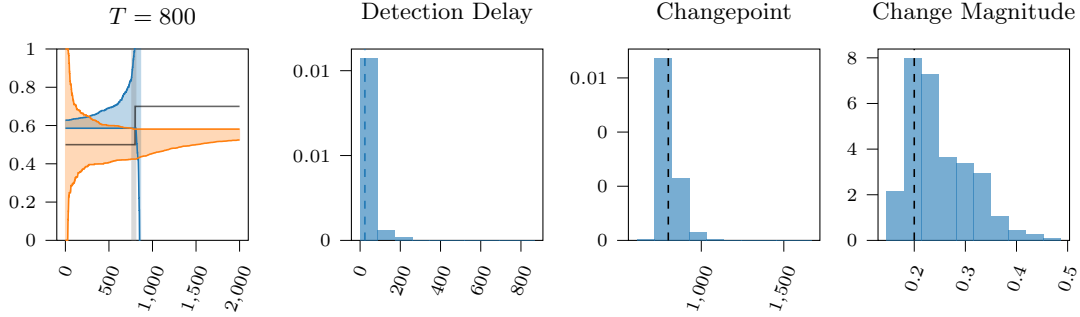} 
            \input{Figures/BernoulliMean/BernoulliMean7277DetectionDelay} 
            \input{Figures/BernoulliMean/BernoulliMean7277EstimatedChangepoint}  
            \input{Figures/BernoulliMean/BernoulliMean7277EstimatedChangeMagnitude}
            \caption{The figures show the performance of our changepoint detection scheme with independent bounded observations whose mean changes from $p_0=0.4$ to $p_1=0.6$ at the time $T=800$. The first plot shows the forward and backward CSs at time of detection~($\tau=863$) in one of the trials, with the shaded gray region being the points at which the two CSs disagree. The next three plots show the empirical distribution of the detection delay, the estimated changepoint location, and the estimated changepoint magnitude over $250$ trials of the experiment with the same value of $\Delta=0.2$ and $\alpha=0.01$. }
            \label{fig:Bernoulli-mean-CS}
    \end{figure}
    
        \begin{figure}[tb]
            \def\figwidth{0.25\linewidth}
            \def\figheight{0.25\linewidth} %
            \centering
            \input{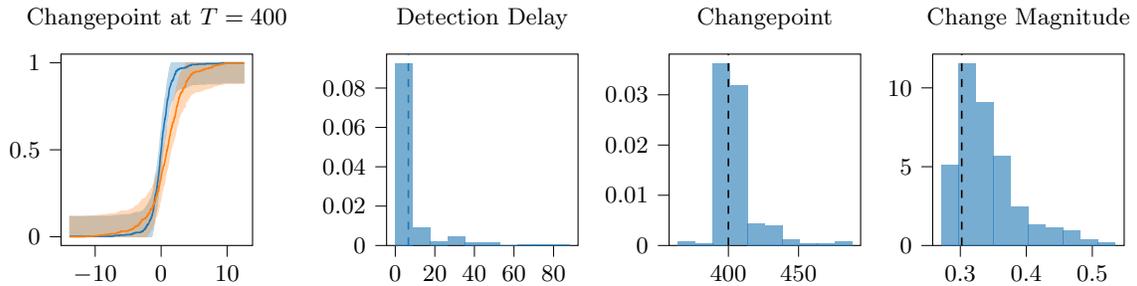} 
            \input{Figures/CDF/CDF3198DetectionDelay} 
            \input{Figures/CDF/CDF3198EstimatedChangepoint}  
            \input{Figures/CDF/CDF3198EstimatedChangeMagnitude}
            \caption{The figures show the performance of our changepoint detection scheme with observations drawn from univariate $t$-distributions~($3$ degrees of freedom) whose mean changes from $0$ to $1.0$ at the time $T=800$. The first plot shows the forward and backward CSs around the empirical CDFs, at time of detection in one of the trials. The next three plots show the empirical distribution of the detection delay, the estimated changepoint location, and the estimated changepoint magnitude over $250$ trials of the experiment.}
            \label{fig:CDF-CS}
    \end{figure}
 
     \begin{figure}
             \def\figwidth{0.25\linewidth}
             \def\figheight{0.25\linewidth} %
            \centering
            \input{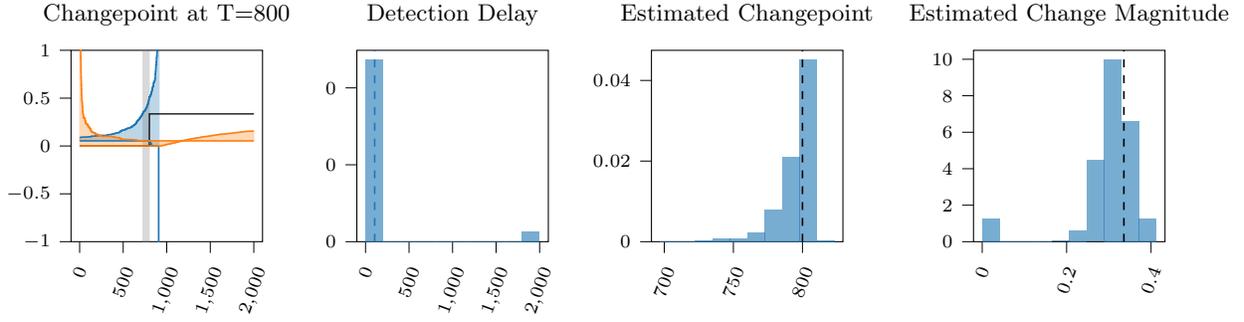} 
            \input{Figures/TwoSample/TwoSample9439DetectionDelay} 
            \input{Figures/TwoSample/TwoSample9439EstimatedChangepoint}
            \input{Figures/TwoSample/TwoSample9439EstimatedChangeMagnitude}
            \caption{The figures show the performance of our changepoint detection scheme with independent paired multivariate-Gaussian observations whose kernel-MMD distance changes from a pre-change value of $0$ to $\Delta \approx 0.33$ at the time $T=800$.  }
            \label{fig:TwoSample-mean-CS}
    \end{figure}
 
     \begin{figure}
         \def\figwidth{0.25\linewidth}
         \def\figheight{0.25\linewidth} %
        \centering
        \input{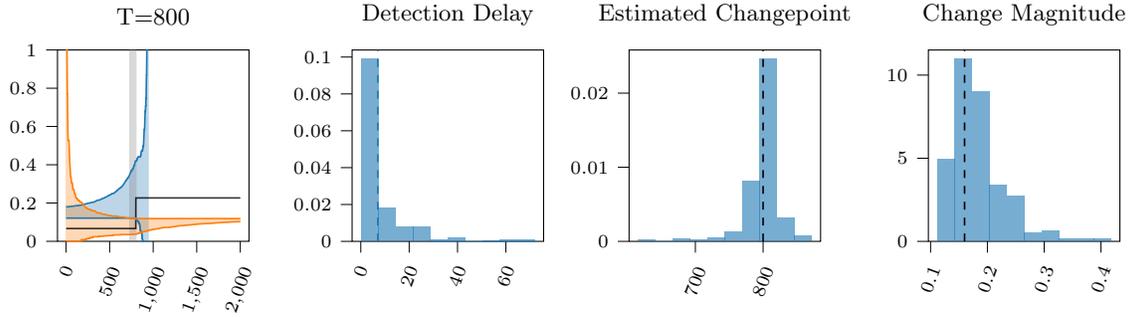} 
        \input{Figures/DistributionShift/DistributionShift3883DetectionDelay} 
        \input{Figures/DistributionShift/DistributionShift3883EstimatedChangepoint}
        \input{Figures/DistributionShift/DistributionShift3883EstimatedChangeMagnitude}
        \caption{The figures show the performance of our scheme for detecting harmful changes in test distribution for two-dimensional feature vectors, as described in~\Cref{subsec:harmful-shifts}. In these plots, there is a change with magnitude $\Delta \approx 0.16$ at the time $T=800$.}
        \label{fig:Distribution-shift-CS}
    \end{figure}
     \begin{figure}
         \def\figwidth{0.32\linewidth}
             \def\figheight{0.32\linewidth} %
            \centering
            \input{Figures/BinaryClassifier/Source} \hspace{1em}
            \input{Figures/BinaryClassifier/Source_new}\hspace{1em}
            \input{Figures/BinaryClassifier/Target}
            \caption{The first plot shows the samples corresponding to the two labels ($L=0$ and $L=1$), as well as the optimal linear classifier $h^*$ for this problem (the dashed black line). In the second plot, the distributions are more separated, and the same classifier $h^*$ is also Bayes-optimal for this problem, with a smaller risk. Finally, in the third figure, we have an example of a harmful distribution shift. In this case, the feature distributions for the two labels are rotated anti-clockwise by $45$ degrees, which makes $h^*$ a suboptimal classifier for this problem. The new Bayes-optimal classifier is shown by the red dotted line.}
            \label{fig:harmful-distribution-shift-1}
    \end{figure}
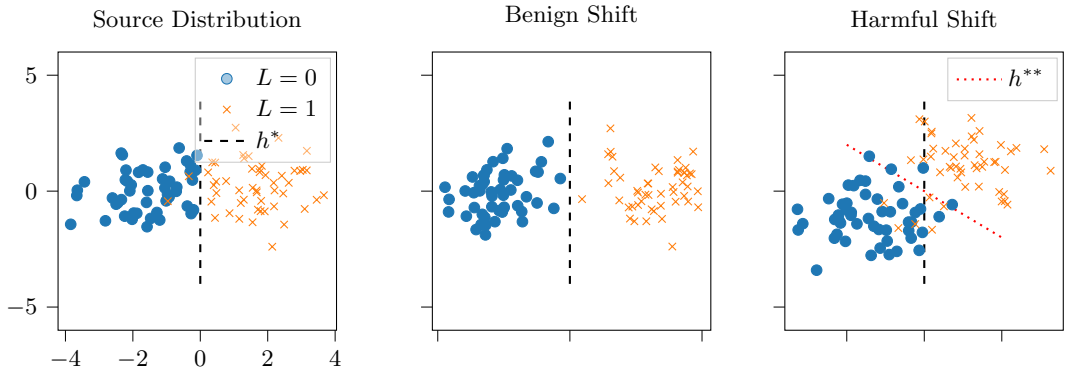
    
    \fi 
    
    \noindent\textbf{Details of the bounded source with a specified mean~(\Cref{subsec:bounded-mean}).} For testing the performance of the SCD scheme described in~\Cref{subsec:bounded-mean}, we constructed a probability distribution, $P_{\theta}$, over $\mc{X} = [0,1]$ with a specified mean $\theta \in [0,1]$ by appropriately mixing two uniform distributions. In particular, we define $P_{\theta} = (1-\theta) U_1 + \theta U_2$, where $U_1 \sim \text{Uniform}([0,\theta])$ and $U_2 \sim \text{Uniform}([\theta,1])$. 
    
    \noindent\textbf{Details of the CDF change detection experiment.~(\Cref{subsec:cdfs}).} For this experiment, we used univariate $t$-distributions with $3$ degrees-of-freedom. For the pre-change distribution, we set the mean to $0$, and the scale parameter to $1$. For the post-change distribution, we set the mean to some value $\delta>0$, and the scale parameter to $2$. 
    
    \noindent\textbf{Details of the Binary classification source~(\Cref{subsec:harmful-shifts})} We consider feature-label pairs $(Z_i, L_i) \in \mathbb{R}^2 \times \{0,1\}$, with a source distribution $P_S = P_L \times P_{Z|L}$. We assume that the label $L$ is drawn uniformly on the set $\{0,1\}$, and the features are drawn from a bivariate normal conditioned on the labels: $P_{Z|L} = N(\mu_L, I_2)$, with $\mu_L = (2L-1) [1, 0]^T \in \mathbb{R}^2$.  For this problem we will consider linear classifiers parameterized by a weight vector $w \in\mathbb{R}^2$, of the form $h_w(z) = \boldsymbol{1}_{ \langle w, z \rangle \geq 0}$. For the $0$-$1$ loss function, $\phi(z, l, h_w) = \boldsymbol{1}_{h_w(z)\neq l}$, it is easy to check that the Bayes-optimal classifier for the source distribution is $h^* \equiv h_{w^*}$, with $w^* \propto [1,0]^T$. For the post change distributions, we rotate the mean of the features by an angle $\gamma$; that is, $\mu_L = (2L-1) [\cos\gamma, \sin \gamma]^T$. 
    \newpage          
    \section{Repeated sequential test interpretation} 
    \label{appendix:repeated-sequential-tests}
        \noindent\textbf{Our scheme as repeated sequential tests.} Due to the equivalence between sequential hypothesis tests and confidence sequences, we can also motivate our general strategy~(\Cref{def:general-strategy}) in the language of sequential hypothesis testing. In particular, our approach can be informally described as follows, due to the time-uniform coverage guarantees of CSs: in each round $t\geq 2$, we run a new sequential hypothesis test for every  $1 \leq s \leq t-1$ to  decide whether $(X_1, \ldots, X_{s})$ and $(X_{s+1}, \ldots, X_t)$ are drawn from the same distribution, or not. As soon as we find a $t$ for which a partition of the observations are sufficiently distinct, we can stop and declare the existence of a changepoint. As we saw in~\Cref{sec:proposed-scheme-2}, this idea can be implemented in an elegant manner by combining a single forward CS~(similar to~\Cref{sec:proposed-scheme-1}) with a succession of CSs constructed on reversed versions of the data, that we called `backward CSs'.  
        
        \noindent\textbf{Connections to CuSum.} We end this section by noting that we can also interpret the popular parametric SCD scheme, CuSum, as also performing repeated sequential tests. 
        Let $f_0$ and $f_1$ denote two density functions on some observation space $\mc{X}$. Let $\{X_t: t \geq 1\}$ denote a sequence of independent observations, and consider a changepoint detection problem with $f_0$ and $f_1$ as the pre- and post-change distributions respectively. 
        
        \begin{definition}[\cusum]
            \label{def:cusum}
            The cumulative sum~(\cusum) method proceeds as follows: 
            \begin{align}
                \tau_c &= \min \{n \geq 1: W_n \geq b_\alpha \}, \quad \text{where} \\
                W_1 &= 0, \quad \text{and} \quad W_n = \max_{1 \leq t \leq n} \prod_{i=t+1}^n \frac{f_1(X_i)}{f_0(X_i)}, \; \text{for } t\geq 2. 
            \end{align}
            The term $b_\alpha$ is selected to ensure that the \arl is at least $1/\alpha$ for a given $\alpha \in (0,1)$. 
        \end{definition}
        
        Observe that the \cusum procedure can also be interpreted as a sequence of repeated sequential~(power-one) tests in the backward direction, testing the null $f_0$ against the alternative $f_1$. These are power-one tests~(as opposed to the usual SPRT) because they only stop when the likelihood ratio process is large, and not when it is small. 
        
        We now introduce an alternative version of \cusum. In this definition, we use the convention that the infimum or minimum over an empty set is infinity, that is,  $\inf \{x: x \in \emptyset\} = \infty$.  
        \begin{definition}[\cusum-II]
            \label{def:modified-cusum}
            For any $n \geq 1$, and $1 \leq t \leq n$, define the backward sequential test, $\backwardtest_n$, as follows: 
           \begin{align}
               \backwardtest_n = \min \{1\leq t \leq n: L_t^n \geq b_\alpha \}, \quad \text{where} \quad L_t^n = \prod_{s=n-t+1}^n \frac{f_1(X_s)}{f_0(X_s)}. 
           \end{align}
           Then, we can define the modified \cusum stopping time as 
           \begin{align}
               \tau_c' \defined \inf \{n: \backwardtest_n<\infty\}. 
           \end{align}
        \end{definition}

        \begin{proposition}
            \label{prop:cusum} The two tests defined in~\Cref{def:cusum} and~\Cref{def:modified-cusum} are the same. 
        \end{proposition}
        
        \begin{proof}
        We show that for any $N \in \mathbb{N}$, the sets $\{\tau_c=N\}$ and $\{\tau_c' = N\}$ are equal. 
        \begin{align}
            \{\tau_c = N\} &= \{W_N\geq b_\alpha\} \cap \lp \cap_{n=1}^{N-1} \{W_n<b_\alpha\} \rp \\
            & = \{\exists t \in [N]: L_t^N \geq b_\alpha \} \cap \lp \cap_{n=1}^{N-1} \{L_t^n <b_\alpha, \; \forall t \in [n]\} \rp \\
            & = \{\backwardtest_N < \infty \} \cap \lp \cap_{n=1}^{N-1}\{\backwardtest_n = \infty \} \rp \\
            & = \{\tau_c' = N\}. 
        \end{align}
        \end{proof}

%% file: Figures/GaussianMean/GaussianMean5564DetectionDelay.tex
\begin{tikzpicture}[baseline]

\definecolor{darkgray176}{RGB}{176,176,176}
\definecolor{steelblue31119180}{RGB}{31,119,180}

\begin{axis}[
height=\figheight,
tick align=outside,
tick pos=left,
title={Detection Delay},
width=\figwidth,
x grid style={darkgray176},
xmin=-100, xmax=2100,
xtick = {0,  1000,  2000},
xtick style={color=black},
yticklabel style={font=\scriptsize},
yticklabels={,,}, 
y grid style={darkgray176},
ymin=0, ymax=0.004452,
ytick style={color=black}
]
\draw[draw=none,fill=steelblue31119180,fill opacity=0.6] (axis cs:0,0) rectangle (axis cs:200,0.00424);
\draw[draw=none,fill=steelblue31119180,fill opacity=0.6] (axis cs:200,0) rectangle (axis cs:400,0.00047);
\draw[draw=none,fill=steelblue31119180,fill opacity=0.6] (axis cs:400,0) rectangle (axis cs:600,0.00011);
\draw[draw=none,fill=steelblue31119180,fill opacity=0.6] (axis cs:600,0) rectangle (axis cs:800,3e-05);
\draw[draw=none,fill=steelblue31119180,fill opacity=0.6] (axis cs:800,0) rectangle (axis cs:1000,1e-05);
\draw[draw=none,fill=steelblue31119180,fill opacity=0.6] (axis cs:1000,0) rectangle (axis cs:1200,0);
\draw[draw=none,fill=steelblue31119180,fill opacity=0.6] (axis cs:1200,0) rectangle (axis cs:1400,0);
\draw[draw=none,fill=steelblue31119180,fill opacity=0.6] (axis cs:1400,0) rectangle (axis cs:1600,0);
\draw[draw=none,fill=steelblue31119180,fill opacity=0.6] (axis cs:1600,0) rectangle (axis cs:1800,0);
\draw[draw=none,fill=steelblue31119180,fill opacity=0.6] (axis cs:1800,0) rectangle (axis cs:2000,0.00014);
\addplot [semithick, steelblue31119180, dashed]
table {%
132.346 -1.0842021724855e-19
132.346 0.004452
};
\end{axis}

\end{tikzpicture}

%% file: Figures/GaussianMean/GaussianMean5564EstimatedChangepoint.tex
\begin{tikzpicture}[baseline]

\definecolor{darkgray176}{RGB}{176,176,176}
\definecolor{steelblue31119180}{RGB}{31,119,180}

\begin{axis}[
height=\figheight,
tick align=outside,
tick pos=left,
title={Changepoint},
width=\figwidth,
x grid style={darkgray176},
xmin=147.4, xmax=1810.6,
xtick = {200,  1000,  1800},
xtick style={color=black},
yticklabel style={font=\scriptsize},
yticklabels={,,}, 
y grid style={darkgray176},
ymin=0, ymax=0.00358653406492913,
ytick style={color=black}
]
\draw[draw=none,fill=steelblue31119180,fill opacity=0.6] (axis cs:223,0) rectangle (axis cs:374.2,1.36085527032029e-05);
\draw[draw=none,fill=steelblue31119180,fill opacity=0.6] (axis cs:374.2,0) rectangle (axis cs:525.4,2.72171054064058e-05);
\draw[draw=none,fill=steelblue31119180,fill opacity=0.6] (axis cs:525.4,0) rectangle (axis cs:676.6,1.36085527032029e-05);
\draw[draw=none,fill=steelblue31119180,fill opacity=0.6] (axis cs:676.6,0) rectangle (axis cs:827.8,0.00341574672850393);
\draw[draw=none,fill=steelblue31119180,fill opacity=0.6] (axis cs:827.8,0) rectangle (axis cs:979,0.00219097698521567);
\draw[draw=none,fill=steelblue31119180,fill opacity=0.6] (axis cs:979,0) rectangle (axis cs:1130.2,0.000639601977050537);
\draw[draw=none,fill=steelblue31119180,fill opacity=0.6] (axis cs:1130.2,0) rectangle (axis cs:1281.4,0.000190519737844841);
\draw[draw=none,fill=steelblue31119180,fill opacity=0.6] (axis cs:1281.4,0) rectangle (axis cs:1432.6,8.16513162192174e-05);
\draw[draw=none,fill=steelblue31119180,fill opacity=0.6] (axis cs:1432.6,0) rectangle (axis cs:1583.8,2.72171054064058e-05);
\draw[draw=none,fill=steelblue31119180,fill opacity=0.6] (axis cs:1583.8,0) rectangle (axis cs:1735,1.36085527032029e-05);
\addplot [semithick, black, dashed]
table {%
800 0
800 0.00358653406492913
};
\end{axis}

\end{tikzpicture}

%% file: Figures/GaussianMean/GaussianMean5564EstimatedChangeMagnitude.tex
\begin{tikzpicture}[baseline]

\definecolor{darkgray176}{RGB}{176,176,176}
\definecolor{steelblue31119180}{RGB}{31,119,180}

\begin{axis}[
height=\figheight,
tick align=outside,
tick pos=left,
title={Change Magnitude},
width=\figwidth,
x grid style={darkgray176},
xmin=-0.105222117999146, xmax=2.20966447798206,
xtick style={color=black},
yticklabel style={font=\scriptsize},
yticklabels={,,},
y grid style={darkgray176},
ymin=0, ymax=2.53463820222821,
ytick style={color=black}
]
\draw[draw=none,fill=steelblue31119180,fill opacity=0.6] (axis cs:0,0) rectangle (axis cs:0.210444235998292,0.133051874132715);
\draw[draw=none,fill=steelblue31119180,fill opacity=0.6] (axis cs:0.210444235998292,0) rectangle (axis cs:0.420888471996583,2.41394114497925);
\draw[draw=none,fill=steelblue31119180,fill opacity=0.6] (axis cs:0.420888471996583,0) rectangle (axis cs:0.631332707994875,1.38754097309831);
\draw[draw=none,fill=steelblue31119180,fill opacity=0.6] (axis cs:0.631332707994875,0) rectangle (axis cs:0.841776943993167,0.304118569446205);
\draw[draw=none,fill=steelblue31119180,fill opacity=0.6] (axis cs:0.841776943993167,0) rectangle (axis cs:1.05222117999146,0.123548168837521);
\draw[draw=none,fill=steelblue31119180,fill opacity=0.6] (axis cs:1.05222117999146,0) rectangle (axis cs:1.26266541598975,0.104540758247133);
\draw[draw=none,fill=steelblue31119180,fill opacity=0.6] (axis cs:1.26266541598975,0) rectangle (axis cs:1.47310965198804,0.0570222317711634);
\draw[draw=none,fill=steelblue31119180,fill opacity=0.6] (axis cs:1.47310965198804,0) rectangle (axis cs:1.68355388798633,0.0665259370663573);
\draw[draw=none,fill=steelblue31119180,fill opacity=0.6] (axis cs:1.68355388798633,0) rectangle (axis cs:1.89399812398462,0.104540758247133);
\draw[draw=none,fill=steelblue31119180,fill opacity=0.6] (axis cs:1.89399812398462,0) rectangle (axis cs:2.10444235998292,0.0570222317711634);
\addplot [semithick, black, dashed]
table {%
0.35 0
0.35 2.53463820222821
};
\end{axis}

\end{tikzpicture}

%% file: Figures/GaussianMean/Delay_vs_Delta8208.tex
\begin{tikzpicture}

\definecolor{darkgray176}{RGB}{176,176,176}
\definecolor{steelblue31119180}{RGB}{31,119,180}
\definecolor{lightgray204}{RGB}{204,204,204}
\definecolor{clr1}{HTML}{FF7F0E} %
\definecolor{clr2}{HTML}{2ca02c} %
\definecolor{clr3}{HTML}{d62728} %
\definecolor{clr4}{HTML}{9467bd} %
\definecolor{clr4}{HTML}{8c564b} %

\begin{axis}[
height=\figheight,
tick align=outside,
tick pos=left,
title={Delay vs Change Magnitude},
width=\figwidth,
x grid style={darkgray176},
xlabel={$\propto 1/\Delta^2$},
xmin=0, xmax=12.45,
xtick style={color=black},
y grid style={darkgray176},
ylabel={Average Detection Delay},
ymin=0.712, ymax=1100,
ytick style={color=black}, 
legend style={fill opacity=0.4, draw opacity=1, text opacity=1, draw=lightgray204, at={(0.03,0.99)},anchor=north west},
]

\addplot [thick, steelblue31119180]
table {%
3 85.58
4 127.5
5 200.12
6 233.28
7 273.4
8 420.7
9 446.8
10 553.86
11 577.82
12 660.5
};
\addlegendentry{Gaussian-mean}; 
\addplot [thick, steelblue31119180, dashed, forget plot]
table {%
3 61.1447272727271
4 127.102787878788
5 193.060848484848
6 259.018909090909
7 324.97696969697
8 390.93503030303
9 456.893090909091
10 522.851151515152
11 588.809212121212
12 654.767272727273
};

\addplot [thick, clr1]
table {%
3 59.96
4 76.24
5 92.92
6 111.4
7 135.64
8 156.96
9 183.6
10 216.16
11 254.28
12 298.5
};
\addlegendentry{Bounded-mean};
\addplot[thick, clr1, dashed, forget plot]
table {%
3 42.73
4 68.47
5 94.21
6 119.95
7 145.69
8 171.43
9 197.17
10 222.91
11 248.65
12 274.39
};
\addplot [thick, clr4]
table {%
1.2 19.1
2.4 34.8
3.6 39.1
4.8 62.6
6 66.4
7.2 106.6
8.4 103.3
9.6 98.5
10.8 141.9
12 139.8
};
\addlegendentry{CDF change}
\addplot [thick, clr4, dashed, forget plot]
table {%
1.2 18.6109090909091
2.4 32.5218181818182
3.6 46.4327272727272
4.8 60.3436363636363
6 74.2545454545454
7.2 88.1654545454545
8.4 102.076363636364
9.6 115.987272727273
10.8 129.898181818182
12 143.809090909091
};
\addplot [thick, clr2]
table {%
4.76622500669398 117.08
5.47057928622595 144.62
6.30403246459839 183.84
6.96922113829242 225.78
8.07751908451832 264.84
8.70376521375936 317.24
9.29447295397706 370.84
10.2199823565354 445.14
10.9865262196266 514.32
12 594.96
};
\addlegendentry{Two-Sample}; 
\addplot [thick, clr2, dashed, forget plot]
table {%
4.76622500669398 84.8721154425327
5.47057928622595 131.587144079766
6.30403246459839 186.864424600206
6.96922113829242 230.981864990329
8.07751908451832 304.487729582977
8.70376521375936 346.022375870483
9.29447295397706 385.200002531375
10.2199823565354 446.582746274014
10.9865262196266 497.422387001115
12 564.639209627203
};

\addplot [thick, clr3]
table {%
12 1058.18
6.61753441518029 792.46
3.8415125236463 394.14
2.33104338111647 255.26
1.47257266106247 155.5
0.966311128539803 100.48
0.657916117370137 75.72
0.464470110592173 50.86
0.33981530570671 44.6
0.257459954641302 35.44
};
\addlegendentry{Distribution-Shift}; 
\addplot [thick, clr3, dashed, forget plot]
table {%
12 1141.92325704574
6.61753441518029 642.015114271581
3.8415125236463 384.186071757854
2.33104338111647 243.897980387726
1.47257266106247 164.1656539055
0.966311128539803 117.145518152996
0.657916117370137 88.5026638987212
0.464470110592173 70.5359472608687
0.33981530570671 58.9583624131518
0.257459954641302 51.3094309058636
};

\end{axis}
\end{tikzpicture}

%% file: Figures/BernoulliMean/BernoulliMean7277DetectionDelay.tex
\begin{tikzpicture}[baseline]

\definecolor{darkgray176}{RGB}{176,176,176}
\definecolor{steelblue31119180}{RGB}{31,119,180}

\begin{axis}[
height=\figheight,
tick align=outside,
tick pos=left,
title={Detection Delay},
width=\figwidth,
x grid style={darkgray176},
xmin=-43.3, xmax=909.3,
xtick style={color=black},
xticklabel style={rotate=70, font=\scriptsize},
yticklabel style={font=\scriptsize},
y grid style={darkgray176},
ymin=0, ymax=0.0112759815242494,
ytick style={color=black}
]
\draw[draw=none,fill=steelblue31119180,fill opacity=0.6] (axis cs:7.105427357601e-15,0) rectangle (axis cs:86.6,0.0107390300230947);
\draw[draw=none,fill=steelblue31119180,fill opacity=0.6] (axis cs:86.6,0) rectangle (axis cs:173.2,0.000577367205542725);
\draw[draw=none,fill=steelblue31119180,fill opacity=0.6] (axis cs:173.2,0) rectangle (axis cs:259.8,0.000184757505773672);
\draw[draw=none,fill=steelblue31119180,fill opacity=0.6] (axis cs:259.8,0) rectangle (axis cs:346.4,0);
\draw[draw=none,fill=steelblue31119180,fill opacity=0.6] (axis cs:346.4,0) rectangle (axis cs:433,2.3094688221709e-05);
\draw[draw=none,fill=steelblue31119180,fill opacity=0.6] (axis cs:433,0) rectangle (axis cs:519.6,0);
\draw[draw=none,fill=steelblue31119180,fill opacity=0.6] (axis cs:519.6,0) rectangle (axis cs:606.2,0);
\draw[draw=none,fill=steelblue31119180,fill opacity=0.6] (axis cs:606.2,0) rectangle (axis cs:692.8,0);
\draw[draw=none,fill=steelblue31119180,fill opacity=0.6] (axis cs:692.8,0) rectangle (axis cs:779.4,0);
\draw[draw=none,fill=steelblue31119180,fill opacity=0.6] (axis cs:779.4,0) rectangle (axis cs:866,2.3094688221709e-05);
\addplot [semithick, steelblue31119180, dashed]
table {%
23.808 0
23.808 0.0112759815242494
};
\end{axis}

\end{tikzpicture}

%% file: Figures/BernoulliMean/BernoulliMean7277EstimatedChangepoint.tex
\begin{tikzpicture}[baseline]

\definecolor{darkgray176}{RGB}{176,176,176}
\definecolor{steelblue31119180}{RGB}{31,119,180}

\begin{axis}[
height=\figheight,
tick align=outside,
tick pos=left,
title={Changepoint},
width=\figwidth,
x grid style={darkgray176},
xmin=558.25, xmax=1718.75,
xtick style={color=black},
xticklabel style={rotate=70, font=\scriptsize},
yticklabel style={font=\scriptsize},
y grid style={darkgray176},
ymin=0, ymax=0.00706635071090047,
ytick style={color=black}
]
\draw[draw=none,fill=steelblue31119180,fill opacity=0.6] (axis cs:611,0) rectangle (axis cs:716.5,5.68720379146919e-05);
\draw[draw=none,fill=steelblue31119180,fill opacity=0.6] (axis cs:716.5,0) rectangle (axis cs:822,0.00672985781990521);
\draw[draw=none,fill=steelblue31119180,fill opacity=0.6] (axis cs:822,0) rectangle (axis cs:927.5,0.00231279620853081);
\draw[draw=none,fill=steelblue31119180,fill opacity=0.6] (axis cs:927.5,0) rectangle (axis cs:1033,0.000303317535545024);
\draw[draw=none,fill=steelblue31119180,fill opacity=0.6] (axis cs:1033,0) rectangle (axis cs:1138.5,3.7914691943128e-05);
\draw[draw=none,fill=steelblue31119180,fill opacity=0.6] (axis cs:1138.5,0) rectangle (axis cs:1244,1.8957345971564e-05);
\draw[draw=none,fill=steelblue31119180,fill opacity=0.6] (axis cs:1244,0) rectangle (axis cs:1349.5,0);
\draw[draw=none,fill=steelblue31119180,fill opacity=0.6] (axis cs:1349.5,0) rectangle (axis cs:1455,0);
\draw[draw=none,fill=steelblue31119180,fill opacity=0.6] (axis cs:1455,0) rectangle (axis cs:1560.5,0);
\draw[draw=none,fill=steelblue31119180,fill opacity=0.6] (axis cs:1560.5,0) rectangle (axis cs:1666,1.8957345971564e-05);
\addplot [semithick, black, dashed]
table {%
800 -2.16840434497101e-19
800 0.00706635071090047
};
\end{axis}

\end{tikzpicture}

%% file: Figures/BernoulliMean/BernoulliMean7277EstimatedChangeMagnitude.tex
\begin{tikzpicture}[baseline]

\definecolor{darkgray176}{RGB}{176,176,176}
\definecolor{steelblue31119180}{RGB}{31,119,180}

\begin{axis}[
height=\figheight,
tick align=outside,
tick pos=left,
title={Change Magnitude},
width=\figwidth,
x grid style={darkgray176},
xmin=0.129147723750587, xmax=0.503748453784166,
xtick style={color=black},
xticklabel style={rotate=70, font=\scriptsize},
yticklabel style={font=\scriptsize},
y grid style={darkgray176},
ymin=0, ymax=8.38652930473036,
ytick style={color=black}
]
\draw[draw=none,fill=steelblue31119180,fill opacity=0.6] (axis cs:0.146175029661204,0) rectangle (axis cs:0.180229641482438,2.17298028203798);
\draw[draw=none,fill=steelblue31119180,fill opacity=0.6] (axis cs:0.180229641482438,0) rectangle (axis cs:0.214284253303673,7.98717076640986);
\draw[draw=none,fill=steelblue31119180,fill opacity=0.6] (axis cs:0.214284253303673,0) rectangle (axis cs:0.248338865124907,7.28242040466782);
\draw[draw=none,fill=steelblue31119180,fill opacity=0.6] (axis cs:0.248338865124907,0) rectangle (axis cs:0.282393476946142,3.6412102023339);
\draw[draw=none,fill=steelblue31119180,fill opacity=0.6] (axis cs:0.282393476946142,0) rectangle (axis cs:0.316448088767376,3.40629341508656);
\draw[draw=none,fill=steelblue31119180,fill opacity=0.6] (axis cs:0.316448088767376,0) rectangle (axis cs:0.350502700588611,2.93645984059186);
\draw[draw=none,fill=steelblue31119180,fill opacity=0.6] (axis cs:0.350502700588611,0) rectangle (axis cs:0.384557312409845,1.05712554261307);
\draw[draw=none,fill=steelblue31119180,fill opacity=0.6] (axis cs:0.384557312409845,0) rectangle (axis cs:0.41861192423108,0.469833574494697);
\draw[draw=none,fill=steelblue31119180,fill opacity=0.6] (axis cs:0.41861192423108,0) rectangle (axis cs:0.452666536052314,0.293645984059186);
\draw[draw=none,fill=steelblue31119180,fill opacity=0.6] (axis cs:0.452666536052314,0) rectangle (axis cs:0.486721147873549,0.117458393623674);
\addplot [semithick, black, dashed]
table {%
0.2 0
0.2 8.38652930473036
};
\end{axis}

\end{tikzpicture}

%% file: Figures/CDF/CDF3198DetectionDelay.tex
\begin{tikzpicture}[baseline]

\definecolor{darkgray176}{RGB}{176,176,176}
\definecolor{steelblue31119180}{RGB}{31,119,180}

\begin{axis}[
height=\figheight,
tick align=outside,
tick pos=left,
title={Detection Delay},
width=\figwidth,
x grid style={darkgray176},
xmin=-4.4, xmax=92.4,
xtick style={color=black},
y grid style={darkgray176},
ymin=0, ymax=0.0972443181818182,
ytick style={color=black}
]
\draw[draw=none,fill=steelblue31119180,fill opacity=0.6] (axis cs:0,0) rectangle (axis cs:8.8,0.0926136363636363);
\draw[draw=none,fill=steelblue31119180,fill opacity=0.6] (axis cs:8.8,0) rectangle (axis cs:17.6,0.00909090909090909);
\draw[draw=none,fill=steelblue31119180,fill opacity=0.6] (axis cs:17.6,0) rectangle (axis cs:26.4,0.00227272727272727);
\draw[draw=none,fill=steelblue31119180,fill opacity=0.6] (axis cs:26.4,0) rectangle (axis cs:35.2,0.00454545454545455);
\draw[draw=none,fill=steelblue31119180,fill opacity=0.6] (axis cs:35.2,0) rectangle (axis cs:44,0.00170454545454545);
\draw[draw=none,fill=steelblue31119180,fill opacity=0.6] (axis cs:44,0) rectangle (axis cs:52.8,0.00170454545454545);
\draw[draw=none,fill=steelblue31119180,fill opacity=0.6] (axis cs:52.8,0) rectangle (axis cs:61.6,0);
\draw[draw=none,fill=steelblue31119180,fill opacity=0.6] (axis cs:61.6,0) rectangle (axis cs:70.4,0.000568181818181818);
\draw[draw=none,fill=steelblue31119180,fill opacity=0.6] (axis cs:70.4,0) rectangle (axis cs:79.2,0.000568181818181818);
\draw[draw=none,fill=steelblue31119180,fill opacity=0.6] (axis cs:79.2,0) rectangle (axis cs:88,0.000568181818181818);
\addplot [semithick, steelblue31119180, dashed]
table {%
6.62 0
6.62 0.0972443181818182
};
\end{axis}

\end{tikzpicture}

%% file: Figures/CDF/CDF3198EstimatedChangepoint.tex
\begin{tikzpicture}[baseline]

\definecolor{darkgray176}{RGB}{176,176,176}
\definecolor{steelblue31119180}{RGB}{31,119,180}

\begin{axis}[
height=\figheight,
tick align=outside,
tick pos=left,
title={Changepoint},
width=\figwidth,
x grid style={darkgray176},
xmin=357.8, xmax=494.2,
xtick style={color=black},
y grid style={darkgray176},
ymin=0, ymax=0.0381048387096775,
ytick style={color=black}
]
\draw[draw=none,fill=steelblue31119180,fill opacity=0.6] (axis cs:364,0) rectangle (axis cs:376.4,0.000806451612903227);
\draw[draw=none,fill=steelblue31119180,fill opacity=0.6] (axis cs:376.4,0) rectangle (axis cs:388.8,0.000403225806451612);
\draw[draw=none,fill=steelblue31119180,fill opacity=0.6] (axis cs:388.8,0) rectangle (axis cs:401.2,0.0362903225806452);
\draw[draw=none,fill=steelblue31119180,fill opacity=0.6] (axis cs:401.2,0) rectangle (axis cs:413.6,0.0318548387096773);
\draw[draw=none,fill=steelblue31119180,fill opacity=0.6] (axis cs:413.6,0) rectangle (axis cs:426,0.00443548387096775);
\draw[draw=none,fill=steelblue31119180,fill opacity=0.6] (axis cs:426,0) rectangle (axis cs:438.4,0.00403225806451614);
\draw[draw=none,fill=steelblue31119180,fill opacity=0.6] (axis cs:438.4,0) rectangle (axis cs:450.8,0.00120967741935484);
\draw[draw=none,fill=steelblue31119180,fill opacity=0.6] (axis cs:450.8,0) rectangle (axis cs:463.2,0.000403225806451614);
\draw[draw=none,fill=steelblue31119180,fill opacity=0.6] (axis cs:463.2,0) rectangle (axis cs:475.6,0.000403225806451612);
\draw[draw=none,fill=steelblue31119180,fill opacity=0.6] (axis cs:475.6,0) rectangle (axis cs:488,0.000806451612903227);
\addplot [semithick, black, dashed]
table {%
400 -8.67361737988404e-19
400 0.0381048387096775
};
\end{axis}

\end{tikzpicture}

%% file: Figures/CDF/CDF3198EstimatedChangeMagnitude.tex
\begin{tikzpicture}[baseline]

\definecolor{darkgray176}{RGB}{176,176,176}
\definecolor{steelblue31119180}{RGB}{31,119,180}

\begin{axis}[
height=\figheight,
tick align=outside,
tick pos=left,
title={Change Magnitude},
width=\figwidth,
x grid style={darkgray176},
xmin=0.258139062893703, xmax=0.548257431768666,
xtick style={color=black},
y grid style={darkgray176},
ymin=0, ymax=12.1424576239716,
ytick style={color=black}
]
\draw[draw=none,fill=steelblue31119180,fill opacity=0.6] (axis cs:0.271326261478929,0) rectangle (axis cs:0.29770065864938,5.11860040354774);
\draw[draw=none,fill=steelblue31119180,fill opacity=0.6] (axis cs:0.29770065864938,0) rectangle (axis cs:0.324075055819831,11.5642453561634);
\draw[draw=none,fill=steelblue31119180,fill opacity=0.6] (axis cs:0.324075055819831,0) rectangle (axis cs:0.350449452990282,9.09973405075154);
\draw[draw=none,fill=steelblue31119180,fill opacity=0.6] (axis cs:0.350449452990282,0) rectangle (axis cs:0.376823850160734,5.68733378171973);
\draw[draw=none,fill=steelblue31119180,fill opacity=0.6] (axis cs:0.376823850160734,0) rectangle (axis cs:0.403198247331185,2.46451130541188);
\draw[draw=none,fill=steelblue31119180,fill opacity=0.6] (axis cs:0.403198247331185,0) rectangle (axis cs:0.429572644501636,1.32704454906793);
\draw[draw=none,fill=steelblue31119180,fill opacity=0.6] (axis cs:0.429572644501636,0) rectangle (axis cs:0.455947041672087,1.13746675634395);
\draw[draw=none,fill=steelblue31119180,fill opacity=0.6] (axis cs:0.455947041672087,0) rectangle (axis cs:0.482321438842538,0.947888963619953);
\draw[draw=none,fill=steelblue31119180,fill opacity=0.6] (axis cs:0.482321438842538,0) rectangle (axis cs:0.50869583601299,0.379155585447981);
\draw[draw=none,fill=steelblue31119180,fill opacity=0.6] (axis cs:0.508695836012989,0) rectangle (axis cs:0.535070233183441,0.189577792723991);
\addplot [semithick, black, dashed]
table {%
0.30214 0
0.30214 12.1424576239716
};
\end{axis}

\end{tikzpicture}

%% file: Figures/TwoSample/TwoSample9439DetectionDelay.tex
\begin{tikzpicture}[baseline]

\definecolor{darkgray176}{RGB}{176,176,176}
\definecolor{steelblue31119180}{RGB}{31,119,180}

\begin{axis}[
height=\figheight,
tick align=outside,
tick pos=left,
title={Detection Delay},
width=\figwidth,
x grid style={darkgray176},
xmin=-100, xmax=2100,
xtick style={color=black},
xticklabel style={rotate=70, font=\scriptsize},
yticklabel style={font=\scriptsize},
y grid style={darkgray176},
ymin=0, ymax=0.004977,
ytick style={color=black}
]
\draw[draw=none,fill=steelblue31119180,fill opacity=0.6] (axis cs:0,0) rectangle (axis cs:200,0.00474);
\draw[draw=none,fill=steelblue31119180,fill opacity=0.6] (axis cs:200,0) rectangle (axis cs:400,0);
\draw[draw=none,fill=steelblue31119180,fill opacity=0.6] (axis cs:400,0) rectangle (axis cs:600,0);
\draw[draw=none,fill=steelblue31119180,fill opacity=0.6] (axis cs:600,0) rectangle (axis cs:800,0);
\draw[draw=none,fill=steelblue31119180,fill opacity=0.6] (axis cs:800,0) rectangle (axis cs:1000,0);
\draw[draw=none,fill=steelblue31119180,fill opacity=0.6] (axis cs:1000,0) rectangle (axis cs:1200,0);
\draw[draw=none,fill=steelblue31119180,fill opacity=0.6] (axis cs:1200,0) rectangle (axis cs:1400,0);
\draw[draw=none,fill=steelblue31119180,fill opacity=0.6] (axis cs:1400,0) rectangle (axis cs:1600,0);
\draw[draw=none,fill=steelblue31119180,fill opacity=0.6] (axis cs:1600,0) rectangle (axis cs:1800,0);
\draw[draw=none,fill=steelblue31119180,fill opacity=0.6] (axis cs:1800,0) rectangle (axis cs:2000,0.00026);
\addplot [semithick, steelblue31119180, dashed]
table {%
105.004 0
105.004 0.004977
};
\end{axis}

\end{tikzpicture}

%% file: Figures/TwoSample/TwoSample9439EstimatedChangepoint.tex
\begin{tikzpicture}[baseline]

\definecolor{darkgray176}{RGB}{176,176,176}
\definecolor{steelblue31119180}{RGB}{31,119,180}

\begin{axis}[
height=\figheight,
tick align=outside,
tick pos=left,
title={Estimated Changepoint},
width=\figwidth,
x grid style={darkgray176},
xmin=690.7, xmax=829.3,
xtick style={color=black},
xticklabel style={rotate=70, font=\scriptsize},
yticklabel style={font=\scriptsize},
y grid style={darkgray176},
ymin=0, ymax=0.0474683544303797,
ytick style={color=black}
]
\draw[draw=none,fill=steelblue31119180,fill opacity=0.6] (axis cs:697,0) rectangle (axis cs:709.6,0.000167436876297636);
\draw[draw=none,fill=steelblue31119180,fill opacity=0.6] (axis cs:709.6,0) rectangle (axis cs:722.2,0.000167436876297636);
\draw[draw=none,fill=steelblue31119180,fill opacity=0.6] (axis cs:722.2,0) rectangle (axis cs:734.8,0.000334873752595274);
\draw[draw=none,fill=steelblue31119180,fill opacity=0.6] (axis cs:734.8,0) rectangle (axis cs:747.4,0.000837184381488177);
\draw[draw=none,fill=steelblue31119180,fill opacity=0.6] (axis cs:747.4,0) rectangle (axis cs:760,0.000837184381488177);
\draw[draw=none,fill=steelblue31119180,fill opacity=0.6] (axis cs:760,0) rectangle (axis cs:772.6,0.0023441162681669);
\draw[draw=none,fill=steelblue31119180,fill opacity=0.6] (axis cs:772.6,0) rectangle (axis cs:785.2,0.0080369700622865);
\draw[draw=none,fill=steelblue31119180,fill opacity=0.6] (axis cs:785.2,0) rectangle (axis cs:797.8,0.0210970464135023);
\draw[draw=none,fill=steelblue31119180,fill opacity=0.6] (axis cs:797.8,0) rectangle (axis cs:810.4,0.0452079566003616);
\draw[draw=none,fill=steelblue31119180,fill opacity=0.6] (axis cs:810.4,0) rectangle (axis cs:823,0.000334873752595271);
\addplot [semithick, black, dashed]
table {%
800 -8.67361737988404e-19
800 0.0474683544303797
};
\end{axis}

\end{tikzpicture}

%% file: Figures/TwoSample/TwoSample9439EstimatedChangeMagnitude.tex
\begin{tikzpicture}[baseline]

\definecolor{darkgray176}{RGB}{176,176,176}
\definecolor{steelblue31119180}{RGB}{31,119,180}

\begin{axis}[
height=\figheight,
tick align=outside,
tick pos=left,
title={Estimated Change Magnitude},
width=\figwidth,
x grid style={darkgray176},
xmin=-0.0206206958939851, xmax=0.433034613773687,
xtick style={color=black},
xticklabel style={rotate=70, font=\scriptsize},
yticklabel style={font=\scriptsize},
y grid style={darkgray176},
ymin=0, ymax=10.4894617093448,
ytick style={color=black}
]
\draw[draw=none,fill=steelblue31119180,fill opacity=0.6] (axis cs:3.46944695195361e-18,0) rectangle (axis cs:0.0412413917879702,1.26086918374001);
\draw[draw=none,fill=steelblue31119180,fill opacity=0.6] (axis cs:0.0412413917879702,0) rectangle (axis cs:0.0824827835759404,0);
\draw[draw=none,fill=steelblue31119180,fill opacity=0.6] (axis cs:0.0824827835759404,0) rectangle (axis cs:0.123724175363911,0);
\draw[draw=none,fill=steelblue31119180,fill opacity=0.6] (axis cs:0.123724175363911,0) rectangle (axis cs:0.164965567151881,0);
\draw[draw=none,fill=steelblue31119180,fill opacity=0.6] (axis cs:0.164965567151881,0) rectangle (axis cs:0.206206958939851,0.0484949686053851);
\draw[draw=none,fill=steelblue31119180,fill opacity=0.6] (axis cs:0.206206958939851,0) rectangle (axis cs:0.247448350727821,0.630434591870006);
\draw[draw=none,fill=steelblue31119180,fill opacity=0.6] (axis cs:0.247448350727821,0) rectangle (axis cs:0.288689742515791,4.46153711169542);
\draw[draw=none,fill=steelblue31119180,fill opacity=0.6] (axis cs:0.288689742515791,0) rectangle (axis cs:0.329931134303762,9.98996353270932);
\draw[draw=none,fill=steelblue31119180,fill opacity=0.6] (axis cs:0.329931134303762,0) rectangle (axis cs:0.371172526091732,6.59531573033237);
\draw[draw=none,fill=steelblue31119180,fill opacity=0.6] (axis cs:0.371172526091732,0) rectangle (axis cs:0.412413917879702,1.26086918374001);
\addplot [semithick, black, dashed]
table {%
0.335629137535567 -2.22044604925031e-16
0.335629137535567 10.4894617093448
};
\end{axis}

\end{tikzpicture}

%% file: Figures/DistributionShift/DistributionShift3883DetectionDelay.tex
\begin{tikzpicture}[baseline]

\definecolor{darkgray176}{RGB}{176,176,176}
\definecolor{steelblue31119180}{RGB}{31,119,180}

\begin{axis}[
height=\figheight,
tick align=outside,
tick pos=left,
title={Detection Delay},
width=\figwidth,
x grid style={darkgray176},
xmin=-3.6, xmax=75.6,
xtick style={color=black},
xticklabel style={rotate=70, font=\scriptsize},
yticklabel style={font=\scriptsize},
y grid style={darkgray176},
ymin=0, ymax=0.103833333333333,
ytick style={color=black}
]
\draw[draw=none,fill=steelblue31119180,fill opacity=0.6] (axis cs:-4.44089209850063e-16,0) rectangle (axis cs:7.2,0.0988888888888889);
\draw[draw=none,fill=steelblue31119180,fill opacity=0.6] (axis cs:7.2,0) rectangle (axis cs:14.4,0.0183333333333333);
\draw[draw=none,fill=steelblue31119180,fill opacity=0.6] (axis cs:14.4,0) rectangle (axis cs:21.6,0.00777777777777778);
\draw[draw=none,fill=steelblue31119180,fill opacity=0.6] (axis cs:21.6,0) rectangle (axis cs:28.8,0.00777777777777778);
\draw[draw=none,fill=steelblue31119180,fill opacity=0.6] (axis cs:28.8,0) rectangle (axis cs:36,0.00111111111111111);
\draw[draw=none,fill=steelblue31119180,fill opacity=0.6] (axis cs:36,0) rectangle (axis cs:43.2,0.00222222222222222);
\draw[draw=none,fill=steelblue31119180,fill opacity=0.6] (axis cs:43.2,0) rectangle (axis cs:50.4,0);
\draw[draw=none,fill=steelblue31119180,fill opacity=0.6] (axis cs:50.4,0) rectangle (axis cs:57.6,0.000555555555555555);
\draw[draw=none,fill=steelblue31119180,fill opacity=0.6] (axis cs:57.6,0) rectangle (axis cs:64.8,0.00111111111111111);
\draw[draw=none,fill=steelblue31119180,fill opacity=0.6] (axis cs:64.8,0) rectangle (axis cs:72,0.00111111111111111);
\addplot [semithick, steelblue31119180, dashed]
table {%
7.012 -1.73472347597681e-18
7.012 0.103833333333333
};
\end{axis}

\end{tikzpicture}

%% file: Figures/DistributionShift/DistributionShift3883EstimatedChangepoint.tex
\begin{tikzpicture}[baseline]

\definecolor{darkgray176}{RGB}{176,176,176}
\definecolor{steelblue31119180}{RGB}{31,119,180}

\begin{axis}[
height=\figheight,
tick align=outside,
tick pos=left,
title={Estimated Changepoint},
width=\figwidth,
x grid style={darkgray176},
xmin=602.15, xmax=884.85,
xtick style={color=black},
xticklabel style={rotate=70, font=\scriptsize},
yticklabel style={font=\scriptsize},
y grid style={darkgray176},
ymin=0, ymax=0.0258210116731517,
ytick style={color=black}
]
\draw[draw=none,fill=steelblue31119180,fill opacity=0.6] (axis cs:615,0) rectangle (axis cs:640.7,0.000155642023346303);
\draw[draw=none,fill=steelblue31119180,fill opacity=0.6] (axis cs:640.7,0) rectangle (axis cs:666.4,0);
\draw[draw=none,fill=steelblue31119180,fill opacity=0.6] (axis cs:666.4,0) rectangle (axis cs:692.1,0.000311284046692606);
\draw[draw=none,fill=steelblue31119180,fill opacity=0.6] (axis cs:692.1,0) rectangle (axis cs:717.8,0.000155642023346304);
\draw[draw=none,fill=steelblue31119180,fill opacity=0.6] (axis cs:717.8,0) rectangle (axis cs:743.5,0.00046692607003891);
\draw[draw=none,fill=steelblue31119180,fill opacity=0.6] (axis cs:743.5,0) rectangle (axis cs:769.2,0.00124513618677043);
\draw[draw=none,fill=steelblue31119180,fill opacity=0.6] (axis cs:769.2,0) rectangle (axis cs:794.9,0.0080933852140078);
\draw[draw=none,fill=steelblue31119180,fill opacity=0.6] (axis cs:794.9,0) rectangle (axis cs:820.6,0.0245914396887159);
\draw[draw=none,fill=steelblue31119180,fill opacity=0.6] (axis cs:820.6,0) rectangle (axis cs:846.3,0.00311284046692608);
\draw[draw=none,fill=steelblue31119180,fill opacity=0.6] (axis cs:846.3,0) rectangle (axis cs:872,0.000778210116731516);
\addplot [semithick, black, dashed]
table {%
800 0
800 0.0258210116731517
};
\end{axis}

\end{tikzpicture}

%% file: Figures/DistributionShift/DistributionShift3883EstimatedChangeMagnitude.tex
\begin{tikzpicture}[baseline]

\definecolor{darkgray176}{RGB}{176,176,176}
\definecolor{steelblue31119180}{RGB}{31,119,180}

\begin{axis}[
height=\figheight,
tick align=outside,
tick pos=left,
title={Change Magnitude},
width=\figwidth,
x grid style={darkgray176},
xmin=0.0961052763328491, xmax=0.433279462939312,
xtick style={color=black},
xticklabel style={rotate=70, font=\scriptsize},
yticklabel style={font=\scriptsize},
y grid style={darkgray176},
ymin=0, ymax=11.509777895689,
ytick style={color=black}
]
\draw[draw=none,fill=steelblue31119180,fill opacity=0.6] (axis cs:0.111431375724052,0) rectangle (axis cs:0.142083574506458,4.95886122490004);
\draw[draw=none,fill=steelblue31119180,fill opacity=0.6] (axis cs:0.142083574506458,0) rectangle (axis cs:0.172735773288863,10.9616932339896);
\draw[draw=none,fill=steelblue31119180,fill opacity=0.6] (axis cs:0.172735773288863,0) rectangle (axis cs:0.203387972071269,9.00424801363429);
\draw[draw=none,fill=steelblue31119180,fill opacity=0.6] (axis cs:0.203387972071269,0) rectangle (axis cs:0.234040170853675,3.39290504861582);
\draw[draw=none,fill=steelblue31119180,fill opacity=0.6] (axis cs:0.234040170853675,0) rectangle (axis cs:0.264692369636081,2.74042330849739);
\draw[draw=none,fill=steelblue31119180,fill opacity=0.6] (axis cs:0.264692369636081,0) rectangle (axis cs:0.295344568418486,0.52198539209474);
\draw[draw=none,fill=steelblue31119180,fill opacity=0.6] (axis cs:0.295344568418486,0) rectangle (axis cs:0.325996767200892,0.652481740118427);
\draw[draw=none,fill=steelblue31119180,fill opacity=0.6] (axis cs:0.325996767200892,0) rectangle (axis cs:0.356648965983298,0.130496348023685);
\draw[draw=none,fill=steelblue31119180,fill opacity=0.6] (axis cs:0.356648965983298,0) rectangle (axis cs:0.387301164765704,0.130496348023685);
\draw[draw=none,fill=steelblue31119180,fill opacity=0.6] (axis cs:0.387301164765704,0) rectangle (axis cs:0.417953363548109,0.130496348023686);
\addplot [semithick, black, dashed]
table {%
0.15982015110801 -2.22044604925031e-16
0.15982015110801 11.509777895689
};
\end{axis}

\end{tikzpicture}

%% file: Figures/BinaryClassifier/Source.tex
\begin{tikzpicture}[baseline]

\definecolor{darkgray176}{RGB}{176,176,176}
\definecolor{darkorange25512714}{RGB}{255,127,14}
\definecolor{lightgray204}{RGB}{204,204,204}
\definecolor{steelblue31119180}{RGB}{31,119,180}

\begin{axis}[
height=\figheight,
legend cell align={left},
legend style={fill opacity=0.4, draw opacity=1, text opacity=1, draw=lightgray204},
title = {Source Distribution},
tick align=outside,
tick pos=left,
width=\figwidth,
x grid style={darkgray176},
xmin=-4.21892021001756, xmax=4.04114309939592,
xtick style={color=black},
y grid style={darkgray176},
ymin=-6, ymax=6,
ytick style={color=black}
]
\addplot [draw=steelblue31119180, fill=steelblue31119180, mark=*, only marks]
table{%
x  y
-3.66345761918229 -0.187946844453284
-1.55575168771501 0.0181528432673259
-1.04263684144072 1.03122848794643
-2.33992828258292 -0.36672133105054
-1.47260036345974 -1.17807566991325
-3.43680413946989 0.401134336687806
-1.88387735929065 -0.956668274734436
-0.232990068801567 0.47907285739484
-1.57912876689682 -1.53492160928436
-0.299530482973982 1.14975378012156
-1.20307209784446 -1.24579935127404
-0.404168156150874 1.30104187976526
-2.20423208097155 0.906224430201079
-1.55980609532597 0.821075397007163
-1.97691288061751 -0.947233311894271
-0.203175196685167 -0.805905202006575
-1.8296304449293 0.812398692545145
-0.271743886336677 -0.970565033137043
-0.925988603435842 -0.25390055120291
-2.44451834446622 -0.495926090268073
-3.8434627868624 -1.43101127213946
-0.215446373389072 -0.841730373360224
-3.64743211849844 0.040243313065593
-2.34552215076615 1.63937940737536
-2.2425629546885 -1.07924281373548
-2.02143909243339 -1.21495172317218
-0.713378382312554 -0.00861412454461097
-0.62711657807649 1.86076123527894
-1.70033340330013 0.914086584563881
-0.298683526777083 0.832897001236509
-2.04898263195167 0.269201156741004
-0.0969835106842001 1.54107492740284
-2.16555696117179 0.428421551234583
-1.11766432344232 0.52543420617098
-1.04904756858637 0.137971447344959
-0.142483521147879 0.899679375667534
-2.22351627770248 0.488813523769986
-0.98398921346699 0.387518716223615
-2.49687279844417 -0.563057014284416
-0.908894679089883 0.426015276746414
-1.28653656037529 -0.921147358909037
-1.01279006816784 -0.419832609540001
-2.81381977927632 -1.27971441659574
-2.09001517886807 0.012738046764001
-0.971721518330363 -0.138202280278549
-2.31213190934105 1.56441603015612
-0.697290530844252 0.172417924362363
-0.385955968892347 -0.644813964610108
-1.59601758666909 -0.485613610875445
-2.60604917123793 -0.296143202545706
};
\addlegendentry{$L=0$}
\addplot [draw=darkorange25512714, fill=darkorange25512714, mark=x, only marks]
table{%
x  y
-0.963174282452546 -0.433426400303064
3.45277350911305 -0.373537651852985
1.93551531343064 0.689798622710234
-0.341545057035229 0.640357993908329
2.46149846552117 -0.101254932307187
2.79392258796013 0.883017849166699
3.66568567624076 -0.176294401247972
2.13464792366848 -2.39687967161019
0.399322889113604 0.258314806680432
0.477913395277075 -0.0825221823122167
1.80535427683087 -0.420604188613493
2.18792567186889 -0.658629569467259
1.20146100476919 -0.51363224045597
2.4173309538236 0.665194361537827
1.75723407922844 0.948011420639851
1.2732765156639 1.55839965358749
1.43503747902314 1.50820506587146
0.436718334200379 -1.15306060448745
1.55404232949058 0.144403978288718
1.825715039928 1.09024629972517
1.66859815445396 0.0315400157518738
3.16176264443211 1.73890945797561
1.91145769683899 -0.874067021635316
0.441927841124807 1.23143857253454
1.70899555654274 -0.798948907764655
3.15761554011252 0.89463821561012
1.67533067865817 -0.0521350784453015
0.232094537624855 0.808647843772596
1.5340561347919 0.965344633276562
1.29981624043916 1.41452250860476
0.329227805730135 -0.453228819771878
0.469339619433925 0.41157365480743
3.03999382234528 0.918502563492872
1.55062613817816 -1.34301233780076
0.355556336686159 1.24304035992173
1.1569743342498 -1.19448196518265
0.867393201278354 -0.951491568853024
2.9641056894565 0.870284643286295
2.48726274059303 -1.44167170757196
1.15347083350094 0.0326381041360011
2.31991004556484 2.29768776214805
0.7665418175693 0.43391941639912
3.09642891001016 -0.768067628799544
1.04824092464276 2.73812165058533
1.29056114180197 -0.369190086587081
2.1771836148702 0.359510859926888
1.78981418931611 -0.862838055919066
1.76107640185309 0.817329722018912
2.62373321722325 0.359336318389052
2.20247446122467 0.0293977531483798
};
\addlegendentry{$L=1$}
\addplot [thick, black, dashed]
table {%
0 -4
0 -3.91919191919192
0 -3.83838383838384
0 -3.75757575757576
0 -3.67676767676768
0 -3.5959595959596
0 -3.51515151515152
0 -3.43434343434343
0 -3.35353535353535
0 -3.27272727272727
0 -3.19191919191919
0 -3.11111111111111
0 -3.03030303030303
0 -2.94949494949495
0 -2.86868686868687
0 -2.78787878787879
0 -2.70707070707071
0 -2.62626262626263
0 -2.54545454545455
0 -2.46464646464646
0 -2.38383838383838
0 -2.3030303030303
0 -2.22222222222222
0 -2.14141414141414
0 -2.06060606060606
0 -1.97979797979798
0 -1.8989898989899
0 -1.81818181818182
0 -1.73737373737374
0 -1.65656565656566
0 -1.57575757575758
0 -1.49494949494949
0 -1.41414141414141
0 -1.33333333333333
0 -1.25252525252525
0 -1.17171717171717
0 -1.09090909090909
0 -1.01010101010101
0 -0.929292929292929
0 -0.848484848484848
0 -0.767676767676767
0 -0.686868686868686
0 -0.606060606060606
0 -0.525252525252525
0 -0.444444444444444
0 -0.363636363636363
0 -0.282828282828282
0 -0.202020202020202
0 -0.121212121212121
0 -0.0404040404040402
0 0.0404040404040407
0 0.121212121212122
0 0.202020202020202
0 0.282828282828283
0 0.363636363636364
0 0.444444444444445
0 0.525252525252526
0 0.606060606060606
0 0.686868686868687
0 0.767676767676768
0 0.848484848484849
0 0.92929292929293
0 1.01010101010101
0 1.09090909090909
0 1.17171717171717
0 1.25252525252525
0 1.33333333333333
0 1.41414141414141
0 1.4949494949495
0 1.57575757575758
0 1.65656565656566
0 1.73737373737374
0 1.81818181818182
0 1.8989898989899
0 1.97979797979798
0 2.06060606060606
0 2.14141414141414
0 2.22222222222222
0 2.3030303030303
0 2.38383838383838
0 2.46464646464647
0 2.54545454545455
0 2.62626262626263
0 2.70707070707071
0 2.78787878787879
0 2.86868686868687
0 2.94949494949495
0 3.03030303030303
0 3.11111111111111
0 3.19191919191919
0 3.27272727272727
0 3.35353535353535
0 3.43434343434344
0 3.51515151515152
0 3.5959595959596
0 3.67676767676768
0 3.75757575757576
0 3.83838383838384
0 3.91919191919192
0 4
};
\addlegendentry{$h^*$}
\end{axis}

\end{tikzpicture}

%% file: Figures/BinaryClassifier/Source_new.tex
\begin{tikzpicture}[baseline]

\definecolor{darkgray176}{RGB}{176,176,176}
\definecolor{darkorange25512714}{RGB}{255,127,14}
\definecolor{lightgray204}{RGB}{204,204,204}
\definecolor{steelblue31119180}{RGB}{31,119,180}

\begin{axis}[
height=\figheight,
legend cell align={left},
legend style={fill opacity=0.8, draw opacity=1, text opacity=1, draw=lightgray204},
title = {Benign Shift},
tick align=outside,
tick pos=left,
width=\figwidth,
x grid style={darkgray176},
xmin=-5.20400187726269, xmax=5.3329086813822,
xtick style={color=black},
xticklabels={,,}, 
yticklabels={,,}, 
y grid style={darkgray176},
ymin=-6, ymax=6,
ytick style={color=black}
]
\addplot [draw=steelblue31119180, fill=steelblue31119180, mark=*, only marks]
table{%
x  y
-2.51714086397113 -0.174775705749558
-3.24365806361913 0.407436212880964
-3.35983342327489 -1.44092478789499
-2.90730542283474 1.26882357004697
-3.2578013719396 0.686462638330615
-2.20303814389862 0.656345084204614
-1.98480999980526 0.469555116713025
-2.97666457810212 -0.0327815243797811
-2.51154607662652 0.485000559917849
-1.64425055708231 -0.243800398156602
-0.60789507750517 -0.744062041516603
-1.20428464234835 0.822404162272153
-1.07988215230897 1.26413870192375
-3.21748317883991 -1.1345183551349
-0.362529925280942 0.542623737464901
-2.48844827513253 -0.695289201079806
-2.63177936059579 0.0665806311425845
-3.42338808707702 0.231864661957778
-2.50085206401009 -0.732147081411646
-1.24188327769664 -0.510398624547746
-3.19476622803942 -1.88726360395836
-3.68461703136155 0.611676233588796
-1.82207604893159 -0.882520850637421
-2.94657329856347 -0.211233669900268
-3.31856823073954 0.274423691305326
-3.14085516059877 0.918554482807844
-3.63647306662034 -0.0688281749494384
-3.54881852874679 -1.65783602986523
-2.23509587034013 0.0492092931352116
-4.59087035138818 -0.896143971337121
-2.47829320727104 0.70362918746324
-2.37413420150741 1.83138038769129
-4.57137513086564 -0.355142906787689
-3.91768830918618 -1.08024692843296
-1.79562273384073 -1.31798988365123
-2.5588837835623 -0.902548081862375
-3.61489510552038 -0.708725273812318
-2.67816081899087 -0.837743023221858
-1.29511009853393 0.854919248945833
-2.53786043265991 1.41654861891911
-3.31406679003044 -1.00076947768468
-4.72505139732429 0.169521195364171
-2.95310850079524 -0.672387752582544
-2.8373926677643 -1.31768933695612
-3.41218356129811 0.01134289509861
-3.96111871562034 0.00690724731382415
-3.25382155724641 0.6670803001568
-0.820139409078144 2.13298198243747
-3.18702643357427 -1.46407797614057
-1.97496892653325 -0.629540425155432
};
\addplot [draw=darkorange25512714, fill=darkorange25512714, mark=x, only marks]
table{%
x  y
4.85395820144379 -0.696381110991995
2.36319071411311 -0.666589731902796
1.86855096333401 -0.400360410501286
4.47595670986647 -0.448699985881822
4.37963562967005 1.00942254916407
3.44575448675362 -0.112116054195521
4.00118104404592 -0.455410716389314
2.57565100611504 -0.0604409373268489
3.88839735430975 0.174968253083674
4.42608468107787 0.680148522233863
4.07545251951786 1.34121488544149
2.8927752116907 -1.30034410592386
1.53451315543887 2.70642027197948
2.44030897811677 -1.30107513112947
3.36486711764533 0.860104687337835
3.1782608878776 -0.859821989254968
2.94265879257936 0.333325543417411
4.54262939625677 0.761994420889911
1.4469843736042 1.59183796676233
4.33661067749269 0.796677762953974
4.07543276716419 0.846066214396952
4.83984980070469 -0.00563590838454376
4.26533719084638 0.379241245398329
2.22448964928413 -0.704033272030417
3.23197257616968 1.02241014135251
4.12985950454143 0.913192213101769
4.66474195473321 1.7057744310338
2.70615774099228 -0.433989317992178
1.67096826122198 0.860626095721506
0.460311908281028 -0.344939459894646
3.48606503979884 -1.19758506686777
3.92335454719196 -0.0378449983151261
3.87614444841709 -2.38988746975234
2.8572260981699 -0.151685137655021
4.33180995339485 0.796984995465315
1.6128648162151 1.48459327248535
4.32387431554519 -0.206037840986033
1.9060990076391 0.581412894066885
3.09323663066059 -0.785722922055885
3.44161861040339 -0.630405093170189
1.54375043955305 1.69776688431436
2.96786055445028 0.32879005430982
4.08121340412764 0.795030976367838
4.45597313248778 0.200581652083783
2.27328977312322 -1.18316000684884
2.54835929490743 0.0125553835209318
4.60731440273176 0.730976958375227
1.45779167871339 0.343191463561187
2.58062893386221 -0.598730226674464
2.97536602423943 -0.167885838212994
};
\addplot [thick, black, dashed]
table {%
0 -4
0 -3.91919191919192
0 -3.83838383838384
0 -3.75757575757576
0 -3.67676767676768
0 -3.5959595959596
0 -3.51515151515152
0 -3.43434343434343
0 -3.35353535353535
0 -3.27272727272727
0 -3.19191919191919
0 -3.11111111111111
0 -3.03030303030303
0 -2.94949494949495
0 -2.86868686868687
0 -2.78787878787879
0 -2.70707070707071
0 -2.62626262626263
0 -2.54545454545455
0 -2.46464646464646
0 -2.38383838383838
0 -2.3030303030303
0 -2.22222222222222
0 -2.14141414141414
0 -2.06060606060606
0 -1.97979797979798
0 -1.8989898989899
0 -1.81818181818182
0 -1.73737373737374
0 -1.65656565656566
0 -1.57575757575758
0 -1.49494949494949
0 -1.41414141414141
0 -1.33333333333333
0 -1.25252525252525
0 -1.17171717171717
0 -1.09090909090909
0 -1.01010101010101
0 -0.929292929292929
0 -0.848484848484848
0 -0.767676767676767
0 -0.686868686868686
0 -0.606060606060606
0 -0.525252525252525
0 -0.444444444444444
0 -0.363636363636363
0 -0.282828282828282
0 -0.202020202020202
0 -0.121212121212121
0 -0.0404040404040402
0 0.0404040404040407
0 0.121212121212122
0 0.202020202020202
0 0.282828282828283
0 0.363636363636364
0 0.444444444444445
0 0.525252525252526
0 0.606060606060606
0 0.686868686868687
0 0.767676767676768
0 0.848484848484849
0 0.92929292929293
0 1.01010101010101
0 1.09090909090909
0 1.17171717171717
0 1.25252525252525
0 1.33333333333333
0 1.41414141414141
0 1.4949494949495
0 1.57575757575758
0 1.65656565656566
0 1.73737373737374
0 1.81818181818182
0 1.8989898989899
0 1.97979797979798
0 2.06060606060606
0 2.14141414141414
0 2.22222222222222
0 2.3030303030303
0 2.38383838383838
0 2.46464646464647
0 2.54545454545455
0 2.62626262626263
0 2.70707070707071
0 2.78787878787879
0 2.86868686868687
0 2.94949494949495
0 3.03030303030303
0 3.11111111111111
0 3.19191919191919
0 3.27272727272727
0 3.35353535353535
0 3.43434343434344
0 3.51515151515152
0 3.5959595959596
0 3.67676767676768
0 3.75757575757576
0 3.83838383838384
0 3.91919191919192
0 4
};
\end{axis}

\end{tikzpicture}

%% file: Figures/BinaryClassifier/Target.tex
\begin{tikzpicture}[baseline]

\definecolor{darkgray176}{RGB}{176,176,176}
\definecolor{darkorange25512714}{RGB}{255,127,14}
\definecolor{lightgray204}{RGB}{204,204,204}
\definecolor{steelblue31119180}{RGB}{31,119,180}

\begin{axis}[
height=\figheight,
legend cell align={left},
legend style={fill opacity=0.8, draw opacity=1, text opacity=1, draw=lightgray204},
title={Harmful Shift},
tick align=outside,
tick pos=left,
width=\figwidth,
x grid style={darkgray176},
xmin=-3.5976095640207, xmax=3.5863811189867,
xtick style={color=black},
yticklabels={,,}, 
xticklabels={,,}, 
y grid style={darkgray176},
ymin=-6, ymax=6,
ytick style={color=black}
]
\addplot [draw=steelblue31119180, fill=steelblue31119180, mark=*, only marks, forget plot]
table{%
x  y
-1.89181994973948 -1.04109272861201
-1.99790096791384 -0.519815027832495
-2.21437780970994 -1.03317000913693
-1.65800212973237 0.465923423151705
-2.16149439310681 -1.37429914223949
-2.46594598764363 -0.324436028007913
-0.328102221724927 -1.33542719159312
-0.0267079128992112 0.998062121539501
-0.405188364541688 -1.33511910483931
-0.0480467521933539 -1.77906700730693
0.383845231398203 -1.10312052904605
-0.65743289460245 0.189273486400524
-2.10944596047355 -0.564407947414908
-3.13449656685064 -1.40477651028564
-0.401340437264751 -1.70406841183333
-2.03353608925944 -2.1658304507171
-3.25330992507681 -1.67396412126891
-1.21109283024616 -0.343157908912279
-1.44065437485021 -1.18035093172983
-1.37111633696432 -2.76981803420377
-1.41868355511248 1.50406201385493
-0.875815423023135 -0.891094631144123
-0.131968530576302 -2.55434799962459
-2.25243978195454 -1.86175771938092
-0.85195811602016 0.946173586308448
-0.337791607005499 -2.02542091233507
-2.31889645074136 -2.02941265313231
-1.73735161171721 -1.4117957022667
-0.904906111547606 -2.73639391900601
-1.17277777012561 -1.65305264158535
-1.1383936612704 -2.4572466844515
-3.27106453297491 -0.783433117212062
-2.30972691433791 -1.22152882286617
-1.3017820688346 -1.51171689121369
-2.01696168378695 -0.805473162748098
-0.956638811498783 -2.15287319046394
-1.08719821597828 -0.882321416318291
-1.51615427171838 -0.142608559312897
-1.8205756202623 0.236641036167311
0.720783839873343 -0.581086962409566
-0.097831342817653 -0.768740533593585
-1.88403236321961 -0.916570867604183
-2.1928424168596 -1.28370430302315
-1.51704474710946 0.432663915754197
-0.707994163076371 -2.59799380056208
-0.561052516763917 -0.541157576350008
-1.95848868909212 0.258978195335675
-0.230908107464087 -0.92150424315057
-2.77755157037164 -3.40944203921186
-0.98041158129431 -1.68279032771851
};

\addplot [draw=darkorange25512714, fill=darkorange25512714, mark=x, only marks, forget plot]
table{%
x  y
1.42605010194405 2.59337237818511
0.505532321214097 -0.658046179307097
0.690050193700596 1.00310380950995
3.25983608794091 0.880163151839709
-0.34453129444009 1.33478715699553
1.80489208713729 1.24850514133448
1.97770906495185 1.99290201362762
1.49433263027537 1.04790152170486
0.183423946277629 2.4897390105416
0.545900439673308 1.737410980886
3.09420191200733 1.80968684846039
0.814311749398126 2.05839193701916
0.135712266922134 -0.271774749661077
1.84052052544407 1.26627164152378
1.22311889115872 2.36217869288875
1.29928497816173 1.1199853040206
0.932580972056892 0.964166006780321
0.818041105242401 1.06431218716721
1.03989980087147 1.5324272976682
1.05723305685243 0.182879644074757
0.761029002486167 1.46587932916687
-1.03843287314898 -0.511710340853281
0.953206290697841 0.807642304948768
2.14645537571096 -0.581959970666351
1.40725344165999 0.195590502918099
-0.233973177650817 -1.41981189271369
0.231412421954351 1.83435559189079
1.04477484219834 1.83449797373894
0.407058084061602 0.81233877550728
0.858235516695164 0.586909067566598
2.07778218629447 -0.433613528940535
1.42295678246874 1.59162217331323
1.27939710956286 2.21712021064188
1.99635833875915 -0.446339174803206
1.35010485818432 0.9825972182866
1.93280890932514 -0.0370945616367109
2.37826175253451 1.86428133039502
2.05439328209432 -0.129816716272156
0.19210773689007 2.56501379222644
0.0123820471596294 3.00211189966435
-0.12610338886955 3.10163435464847
-0.0219009836516513 1.36828533445628
0.199231856329586 -1.66460662552933
2.23440079866546 1.71504942629987
0.169988253998846 1.48914340634394
0.654919896568331 1.05020786697154
2.08173760922952 0.912606418664631
1.54696175955434 1.26750886521053
-0.668393727213283 -1.59206197279689
1.21699357999328 3.15827931955217
};
\addplot [thick, black, dashed, forget plot]
table {%
0 -4
0 -3.91919191919192
0 -3.83838383838384
0 -3.75757575757576
0 -3.67676767676768
0 -3.5959595959596
0 -3.51515151515152
0 -3.43434343434343
0 -3.35353535353535
0 -3.27272727272727
0 -3.19191919191919
0 -3.11111111111111
0 -3.03030303030303
0 -2.94949494949495
0 -2.86868686868687
0 -2.78787878787879
0 -2.70707070707071
0 -2.62626262626263
0 -2.54545454545455
0 -2.46464646464646
0 -2.38383838383838
0 -2.3030303030303
0 -2.22222222222222
0 -2.14141414141414
0 -2.06060606060606
0 -1.97979797979798
0 -1.8989898989899
0 -1.81818181818182
0 -1.73737373737374
0 -1.65656565656566
0 -1.57575757575758
0 -1.49494949494949
0 -1.41414141414141
0 -1.33333333333333
0 -1.25252525252525
0 -1.17171717171717
0 -1.09090909090909
0 -1.01010101010101
0 -0.929292929292929
0 -0.848484848484848
0 -0.767676767676767
0 -0.686868686868686
0 -0.606060606060606
0 -0.525252525252525
0 -0.444444444444444
0 -0.363636363636363
0 -0.282828282828282
0 -0.202020202020202
0 -0.121212121212121
0 -0.0404040404040402
0 0.0404040404040407
0 0.121212121212122
0 0.202020202020202
0 0.282828282828283
0 0.363636363636364
0 0.444444444444445
0 0.525252525252526
0 0.606060606060606
0 0.686868686868687
0 0.767676767676768
0 0.848484848484849
0 0.92929292929293
0 1.01010101010101
0 1.09090909090909
0 1.17171717171717
0 1.25252525252525
0 1.33333333333333
0 1.41414141414141
0 1.4949494949495
0 1.57575757575758
0 1.65656565656566
0 1.73737373737374
0 1.81818181818182
0 1.8989898989899
0 1.97979797979798
0 2.06060606060606
0 2.14141414141414
0 2.22222222222222
0 2.3030303030303
0 2.38383838383838
0 2.46464646464647
0 2.54545454545455
0 2.62626262626263
0 2.70707070707071
0 2.78787878787879
0 2.86868686868687
0 2.94949494949495
0 3.03030303030303
0 3.11111111111111
0 3.19191919191919
0 3.27272727272727
0 3.35353535353535
0 3.43434343434344
0 3.51515151515152
0 3.5959595959596
0 3.67676767676768
0 3.75757575757576
0 3.83838383838384
0 3.91919191919192
0 4
};
\addplot [thick, red, dotted]
table {%
-2 2
-1.95959595959596 1.95959595959596
-1.91919191919192 1.91919191919192
-1.87878787878788 1.87878787878788
-1.83838383838384 1.83838383838384
-1.7979797979798 1.7979797979798
-1.75757575757576 1.75757575757576
-1.71717171717172 1.71717171717172
-1.67676767676768 1.67676767676768
-1.63636363636364 1.63636363636364
-1.5959595959596 1.5959595959596
-1.55555555555556 1.55555555555556
-1.51515151515152 1.51515151515152
-1.47474747474747 1.47474747474747
-1.43434343434343 1.43434343434343
-1.39393939393939 1.39393939393939
-1.35353535353535 1.35353535353535
-1.31313131313131 1.31313131313131
-1.27272727272727 1.27272727272727
-1.23232323232323 1.23232323232323
-1.19191919191919 1.19191919191919
-1.15151515151515 1.15151515151515
-1.11111111111111 1.11111111111111
-1.07070707070707 1.07070707070707
-1.03030303030303 1.03030303030303
-0.98989898989899 0.98989898989899
-0.949494949494949 0.949494949494949
-0.909090909090909 0.909090909090909
-0.868686868686869 0.868686868686869
-0.828282828282828 0.828282828282828
-0.787878787878788 0.787878787878788
-0.747474747474747 0.747474747474747
-0.707070707070707 0.707070707070707
-0.666666666666667 0.666666666666667
-0.626262626262626 0.626262626262626
-0.585858585858586 0.585858585858586
-0.545454545454545 0.545454545454545
-0.505050505050505 0.505050505050505
-0.464646464646465 0.464646464646465
-0.424242424242424 0.424242424242424
-0.383838383838384 0.383838383838384
-0.343434343434343 0.343434343434343
-0.303030303030303 0.303030303030303
-0.262626262626263 0.262626262626263
-0.222222222222222 0.222222222222222
-0.181818181818182 0.181818181818182
-0.141414141414141 0.141414141414141
-0.101010101010101 0.101010101010101
-0.0606060606060606 0.0606060606060606
-0.0202020202020201 0.0202020202020201
0.0202020202020203 -0.0202020202020203
0.060606060606061 -0.060606060606061
0.101010101010101 -0.101010101010101
0.141414141414141 -0.141414141414141
0.181818181818182 -0.181818181818182
0.222222222222222 -0.222222222222222
0.262626262626263 -0.262626262626263
0.303030303030303 -0.303030303030303
0.343434343434343 -0.343434343434343
0.383838383838384 -0.383838383838384
0.424242424242424 -0.424242424242424
0.464646464646465 -0.464646464646465
0.505050505050505 -0.505050505050505
0.545454545454546 -0.545454545454546
0.585858585858586 -0.585858585858586
0.626262626262626 -0.626262626262626
0.666666666666667 -0.666666666666667
0.707070707070707 -0.707070707070707
0.747474747474748 -0.747474747474748
0.787878787878788 -0.787878787878788
0.828282828282829 -0.828282828282829
0.868686868686869 -0.868686868686869
0.909090909090909 -0.909090909090909
0.94949494949495 -0.94949494949495
0.98989898989899 -0.98989898989899
1.03030303030303 -1.03030303030303
1.07070707070707 -1.07070707070707
1.11111111111111 -1.11111111111111
1.15151515151515 -1.15151515151515
1.19191919191919 -1.19191919191919
1.23232323232323 -1.23232323232323
1.27272727272727 -1.27272727272727
1.31313131313131 -1.31313131313131
1.35353535353535 -1.35353535353535
1.39393939393939 -1.39393939393939
1.43434343434343 -1.43434343434343
1.47474747474747 -1.47474747474747
1.51515151515152 -1.51515151515152
1.55555555555556 -1.55555555555556
1.5959595959596 -1.5959595959596
1.63636363636364 -1.63636363636364
1.67676767676768 -1.67676767676768
1.71717171717172 -1.71717171717172
1.75757575757576 -1.75757575757576
1.7979797979798 -1.7979797979798
1.83838383838384 -1.83838383838384
1.87878787878788 -1.87878787878788
1.91919191919192 -1.91919191919192
1.95959595959596 -1.95959595959596
2 -2
};
\addlegendentry{$h^{**}$}; 
\end{axis}
\end{tikzpicture}